\documentclass[final]{siamltex}
\usepackage{graphicx}
\usepackage{subfigure} 
\usepackage{color}

\usepackage{psfrag} 
\usepackage{amsmath} 
\usepackage{amssymb}  
\usepackage{amsfonts}

\newcommand{\R}{\mathbb{R}}  
\newcommand{\C}{\mathbb{C}}

\newcommand{\lbr}{\lbrack}
\newcommand{\rbr}{\rbrack}
\newcommand{\mH}{\mathrm{H}}
\newcommand{\mL}{\mathrm{L}}

\newcommand{\bfx}{\boldsymbol{x}}
\newcommand{\bfy}{\boldsymbol{y}}
\newcommand{\bfn}{\boldsymbol{n}}

\newtheorem{remark}{Remark}

\DeclareMathOperator{\Id}{I}

\begin{document}

\title{An introduction to Multitrace Formulations and Associated
  Domain Decomposition Solvers}

\author{
X. Claeys\footnote[1]{Sorbonne Universit\'es, UPMC Univ Paris 06, CNRS, INRIA,  UMR 7598, Labo. Jacques-Louis Lions, 
\'Equipe Alpines, 4, place Jussieu 75005, Paris, France, claeys@ann.jussieu.fr},
V. Dolean\footnote[2]{Univ. of Strathclyde, Dept. of Mathematics and
  Statistics, Glasgow, UK.  Victorita.Dolean@strath.ac.uk},
M.J. Gander\footnote[3]{Section de Math\'ematiques, Universit\'e de
Gen\`eve, CP 64, 1211 Gen\`eve, Martin.Gander@math.unige.ch} }

\maketitle

\begin{abstract}
{\color{black} Multitrace formulations (MTFs)} are based on a
decomposition of the problem domain into subdomains, and thus domain
decomposition solvers are of interest. The fully rigorous mathematical
MTF can however be daunting for the non-specialist. We introduce in
this paper MTFs on a simple model problem using concepts familiar to
researchers in domain decomposition. This allows us to get a new
understanding of MTFs and a natural block Jacobi iteration, for which
we determine optimal relaxation parameters. We then show how iterative
multitrace formulation solvers are related to a well known domain
decomposition method called optimal Schwarz method: a method which
used Dirichlet to Neumann maps in the transmission condition. We
finally show that the insight gained from the simple model problem
leads to remarkable identities for Calder\'on projectors and related
operators, and the convergence results and optimal choice of the
relaxation parameter we obtained is independent of the geometry, the
space dimension of the problem{\color{black}, and the precise form of
  the spatial elliptic operator, like for optimal Schwarz methods.} We
illustrate our analysis with numerical experiments.
\end{abstract}  
\begin{keywords}
Multitrace formulations, Calder\'on projectors, Dirichlet to Neumann
operators, optimal Schwarz methods
\end{keywords}
\begin{AMS}
  65M55, 65F10, 65N22
\end{AMS}
\pagestyle{myheadings} 
\thispagestyle{plain} 

%\tableofcontents

\section{Introduction}

Multitrace formulations (MTF) for boundary integral equations (BIE)
were developed over the last few years, see
\cite{Hiptmair:2012:MTB,Claeys:2012:ESC, Claeys:2013:MTB}, for the
simulation of electromagnetic transmission problems in piecewise
constant media, and also \cite{Claeys:2013:ASK} for associated
boundary integral methods. MTFs are naturally adapted to the
development of new block preconditioners, as indicated in
\cite{Hiptmair:2013:DDB}, but very little is known so far about such
associated iterative domain decomposition solvers. The first goal of
our presentation (see Section \ref{1dSec}) is to give an elementary
introduction to MTFs, and the associated concepts of representation
formulas and Calder\'on projectors for a simple model problem in one
spatial dimension, in order to make these concepts accessible for
people working in domain decomposition. This approach allows us to get
a complete understanding of the performance of a block Jacobi
iteration for the MTF applied to our model problem, and to determine
the influence of the relaxation parameter on the convergence of the
block Jacobi method. Based on these results, we establish in Section
\ref{MTFOptSchwarzSec} an interesting connection between MTFs with a
well studied class of domain decomposition methods called optimal
Schwarz methods, see
\cite{Nataf:1994:OIC,Gander:1999:OCO,gander2003optimal,gander2011optimal,gander2012optimal},
and \cite{gander2006optimized} for an overview with further
references. Optimal Schwarz methods use Dirichlet to Neumann operators
in their transmission conditions, and are very much related to the
very recent class of sweeping preconditioners,
\cite{engquist2011sweepingPML,engquist2011sweeping}, see also the
earlier work on the same method known under the name AILU (Analytic
Incomplete LU factorization) in \cite{Gander00a,gander2005incomplete},
or frequency filtering \cite{wagner1997adaptive}. To find the
connection with MTFs, we need to generalize first the MTF to the case
of bounded domains and give a formulation of the Calder\'on projectors
in terms of the Dirichlet to Neumann and Neumann to Dirichlet
operators. We then show in Section \ref{MTFGeneral} that the insight
gained for the one dimensional problem holds in much more general
situations. It allows us to discover remarkable properties of the
Calder\'on projector and related operators in higher spatial
dimensions and on various geometries, and the performance of the block
Jacobi iteration and the dependence on the relaxation parameter remain
as we discovered for the one dimensional model problem. We illustrate
our results with numerical experiments that confirm our analysis.

\section{Multitrace Formulations for a Simple 1D Model Problem}\label{1dSec}

Because multitrace formulations need substantial knowledge in
functional analysis, representation formulas and Calder\'on
projectors, we start in this section by explaining these concepts for
a simple model problem in one spatial dimension, without dwelling on
functional analysis issues. The more general formulation and
associated functional analysis framework will be introduced in Section
\ref{MTFGeneral}.

\subsection{Representation Formulas in 1D}\label{ReprForm1D}

We start by examining for some given constant $a>0$ solutions of
the differential equation
\begin{equation}\label{ModelEq1}
  \left\{\begin{array}{rcl}
    \displaystyle-\frac{d^2u}{dx^2} + a^2u&=&0,\quad
    \textrm{in}\; \mathbb{R}\setminus \{0\},\\
    \displaystyle\lim_{\vert x\vert\rightarrow\infty}|u(x)|&=&0.
  \end{array}\right.
\end{equation}
Since the domain on which the equation holds excludes the point $x=0$,
there are non-zero solutions, and to select a particular one, two more
conditions must be imposed on the solution at $x=0$. In transmission
problems and multitrace formulations, one works with solutions that
can be discontinuous, and we thus introduce the notation of jumps
(with convention that the orientation is from $\mathbb{R}_+$ to
$\mathbb{R}_-$),
\begin{equation}
  \begin{array}{rcl}
    \displaystyle[u] & := & u(0_+) -u(0_-),\\
    \displaystyle\left[\frac{du}{dx}\right] & := & 
      \displaystyle-\frac{du}{dx}(0_+) + \frac{du}{dx}(0_-).
  \end{array}
\end{equation}
Imposing both jumps at $x=0$ selects a unique solution of
(\ref{ModelEq1}); solving for example the case where the solution is
continuous, but has a jump of size $\beta$ in the derivative, 
\begin{equation}
  \left\{\begin{array}{rcl}
    \displaystyle -\frac{d^2u}{dx^2} + a^2u&=&0,\quad
      \textrm{in}\; \mathbb{R} \setminus \{ 0\},\\
    \displaystyle [u] &=&0,\\
    \displaystyle \left[\frac{du}{dx}\right] & = &\beta,\\
    \displaystyle\lim_{x\rightarrow\infty}|u(x)|&=&0,
  \end{array}\right.
\end{equation}
we obtain as solution decaying exponentials in each part of the domain
and conditions determining the constants, 
$$
  \begin{array}{rcl}
    u(x) &=& c_+ e^{-ax} \mathbf{1}_{\mathbb{R}_+}+c_- e^{ax}
      \mathbf{1}_{\mathbb{R}_-},\\
    c_+ - c_-  &=& 0,\\
    a(c_++c_-) &=& \beta.
  \end{array}
$$
Solving the linear system for the constants, we find $c_\pm =
\beta/(2a)$ and hence our solution can be written in compact form 
with the so called {\em Green's function $\mathcal{G}$} as 
\begin{equation}\label{ExprGreenKer1D}
  u(x) = \beta\,\mathcal{G}(x),\quad\quad \mathcal{G}(x):= 
    \frac{e^{-a|x|}}{2a}.  
\end{equation}
If we impose instead a jump $\alpha$ in the solution, but continuous
derivatives,
\begin{equation}
  \left\{\begin{array}{rcl}
    \displaystyle -\frac{d^2u}{dx^2} + a^2u&=&0,\, \mathbb{R} \setminus \{ 0\},\\
    \displaystyle [u] &=&\alpha,\\
    \displaystyle \left[\frac{du}{dx}\right] & = &0,\\
    \displaystyle\lim_{x\rightarrow\infty}|u(x)|&=&0.
  \end{array}\right.
\end{equation}
we find by similar calculations
\begin{equation}
  u(x) = \frac{\alpha}{2}\text{sign}(x) e^{-a|x|} 
       = -\alpha \frac{d \mathcal{G}}{dx}(x),
\end{equation}
where $\mathcal{G}$ is again the Green's function from
(\ref{ExprGreenKer1D}).  By linearity, any function $u(x)$ solution to
(\ref{ModelEq1}) with Dirichlet jump $\lbr u\rbr = \alpha$ and
Neumann jump $\lbr du/dx\rbr =\beta$ is thus given by the formula
\begin{equation}\label{eq:repres0}
  \boxed{u(x)=\left[\frac{du}{dx}\right]{\cal G}(x)-[u]\frac{d {\cal G}(x)}{dx},
    \quad\quad\forall x\in \mathbb{R}\setminus\{0\}.}
\end{equation}
This formula is called {\em representation formula} for the
solution.

\subsection{Calder\'on Projectors in 1D}\label{Calderon1d}

If $u$ is any function satisfying $-d^{2}u/dx^{2} + a^{2}u = 0$ on
$\R_{+}$ and $u(x)\to 0$ for $x\to \infty$, then we can extend the
function by zero on the negative real axis, $u(x) = 0$ for $x<0$, and
it then satisfies (\ref{ModelEq1}). Hence the representation formula
(\ref{eq:repres0}) yields $u(x) = (G_{+}\circ
T_{+}(u))(x)$ for $x\in \R_{+}$, where
\begin{equation}\label{eq:tracep}
  G_+\left(\begin{array}{c}\alpha \\
    \beta\end{array}\right):=-\alpha\frac{d {\cal G}(x)}{dx} 
    +\beta\,{\cal G}(x)\;,\quad\quad
    T_+(u):=\left(\begin{array}{c} u(0_+) \\
    -\frac{du}{dx}(0_+) \end{array}\right). 
\end{equation}
Observe that the composition in the reverse order, $T_{+}\circ G_{+}$,
is here a simple $2\times 2$ matrix whose coefficients can be
explicitly computed from the expression of the Green's function
$\mathcal{G}$ given in (\ref{ExprGreenKer1D}).  This yields
\begin{equation}\label{eq:calp}
  \mathbb{P}_+:= T_+ \circ G_+ = \frac{1}{2}
    \Big( \Id + \underbrace{\left[\begin{array}{cc} 0 & 1/a \\ a &
    0\end{array} \right]}_A\Big),\quad\quad A^2  =\Id\quad \Longrightarrow\quad
    \mathbb{P}_+^2 =  \mathbb{P}_+.
\end{equation}
We therefore see that $\mathbb{P}_+$ is a projector, which is called
{\it Calder\'on projector} associated with $\R_{+}$. 

Similarly, if $u$ is any function satisfying $-d^{2}u/dx^{2} +
a^{2}u = 0$ on $\R_{-}$ and $u(x)\to 0$ for $x\to -\infty$ then,
setting $u(x) = 0$ for $x>0$, the representation formula
(\ref{eq:repres0}) can be applied which yields $u(\bfx) = (G_{-}\circ
T_{-}(u))(\bfx)$ for $\bfx\in \R_{-}$ with
\begin{equation}\label{eq:tracem}
  G_-\left(\begin{array}{c}\alpha \\
    \beta\end{array}\right):=\alpha\frac{d {\cal G}(x)}{dx} 
     +\beta\,{\cal G}(x)\;,\quad\quad
   T_-(u):=\left(\begin{array}{c} u(0_-) \\
    \frac{du}{dx}(0_-) \end{array}\right). 
\end{equation}
Computing $\mathbb{P}_{-} := T_{-}\circ G_{-}$ in the same manner as
above, we find that $\mathbb{P}_{-}=(\Id+A)/2$ with the same matrix
$A$ as in (\ref{eq:calp}), and $\mathbb{P}_{-}^{2} =
\mathbb{P}_{-}$ is the {\it Calder\'on projector} associated with
$\R_{-}$.
\begin{remark}\label{CalderonProjectorRemark}
We see that the Calder\'on projector performs a very simple operation:
it takes two arbitrary jumps along the interface, solves the coupled
transmission problem with these jumps, and then returns the Dirichlet
and Neumann trace of the domain the Calder\'on projector is associated
with.
\end{remark}

\subsection{Multitrace Formulation with 2 Subdomains in 1D}\label{2Dom1D}

Suppose now we have a decomposition of $\mathbb{R}$ into two
subdomains $\Omega_1=\R_{-}$ and $\Omega_2=\R_{+}$. Let $T_{1,2}$ be
the trace operators as defined in \eqref{eq:tracep} and
\eqref{eq:tracem} ($T_1=T_{-}$ and $T_2=T_+$) for the subdomains
$\Omega_{1,2}$, and let $\mathbb{P}_{1,2}$ be the corresponding
Calder\'on projectors as defined in \eqref{eq:calp}
($\mathbb{P}_1=\mathbb{P}_-=:\mathbb{P}$ and
$\mathbb{P}_2=\mathbb{P}_+=\mathbb{P}$). Suppose we want to solve the
transmission problem
\begin{equation}\label{ModelEq2}
  \left\{\begin{array}{l}
    \displaystyle -\frac{d^2u}{dx^2} + a^2u=0,\quad
      \textrm{in}\; \mathbb{R} \setminus \{ 0\},\\
    \displaystyle [u] =\alpha, \quad \left[\frac{du}{dx}\right]  = \beta,\\
    \displaystyle\lim_{\vert x\vert\rightarrow-\infty}u(x)=0.
  \end{array}\right.
\end{equation}
The {\it multitrace formulation} introduced in
\cite{Hiptmair:2012:MTB}, which we present in the form with relaxation
parameters from \cite{Hiptmair:2013:DDB} states that $u$ is solution
to (\ref{ModelEq2}) if its traces $U_{1,2}:=(T_{i} u)_{i=1,2}$ verify
the relations
\begin{equation}\label{eq:multi2sd1}
  \left\{\begin{array}{l}
    (\Id-\mathbb{P}_1) U_1
      +\sigma_1  \left(U_1-X U_2 \right)=F_{1},\\[5pt]
    (\Id-\mathbb{P}_2) U_2
      +\sigma_2  \left(U_2-XU_1\right)=F_{2},
  \end{array}\right.
\end{equation}
where $F_{1} = \sigma_{1}\cdot (-\alpha,\beta)^T$, $F_{2} =
\sigma_{2}\cdot (\alpha,\beta)^T$, $\sigma_{1},\sigma_{2}\in\C$ are
some relaxation parameters, and
\begin{equation}\label{eq:xp}
  X:=\left(\begin{array}{cc}
    1 & 0\\
    0 & -1
  \end{array}\right).
\end{equation}
We see that (\ref{eq:multi2sd1}) clearly holds for the solution $u$:
first $(\Id-\mathbb{P}_j) U_j=0$ by Remark
\ref{CalderonProjectorRemark}, since applying the Calder\'on projector
to a solution just gives the solution itself. Second, the relaxation
term on the left gives precisely the jumps weighted by the relaxation,
which we also find in the right hand side functions $F_j$ which
contain as data the jumps $\alpha$, $\beta$ of the transmission
problem. Note also that this formulation would not make sense for
vanishing relaxation parameters, $\sigma_j=0$, $j=1,2$, since then the
jump data $\alpha$, $\beta$ disappear from the problem formulation.

\begin{figure}
  \centering
  \includegraphics[width=0.4\textwidth]{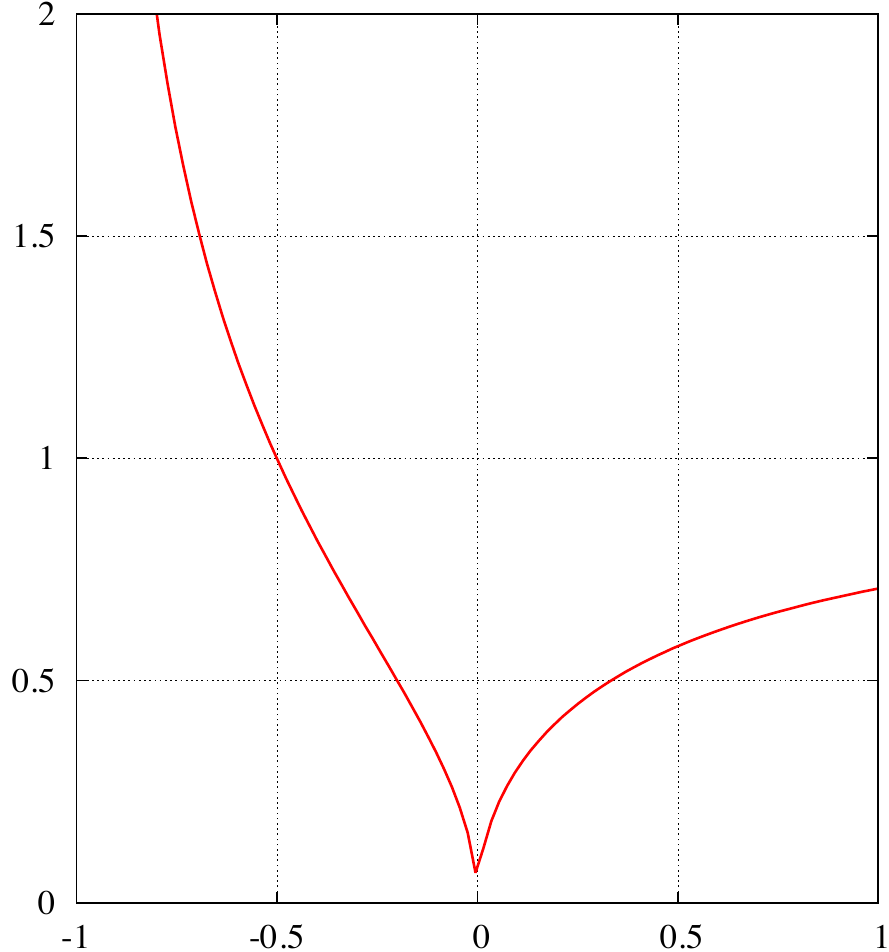}
  \caption{Eigenvalues in modulus of the block Jacobi iteration matrix
    as a function of the relaxation parameter $\sigma_j$}
  \label{fig:conv1d}
\end{figure}

Collecting the operators that act on the same trace variables $U_j$,
we can rewrite (\ref{eq:multi2sd1}) in matrix form as a 
$4\times 4$ linear system of equations, namely
\begin{equation}\label{eq:multi2sd2}
  \left[\begin{array}{cc}
    (1+\sigma_1)\Id-\mathbb{P}_1 & -\sigma_1 X \\
    -\sigma_2 X & (1+\sigma_2)\Id-\mathbb{P}_2
  \end{array}\right] \left[\begin{array}{c} 
    U_1 \\ 
    U_2
  \end{array}\right] = 
  \left[\begin{array}{c} 
    F_{1} \\ 
    F_{2}
  \end{array}\right].
\end{equation}
A very natural iterative method to solve this linear system would
be a block-Jacobi iteration,
\begin{equation}\label{eq:blockJac}
  \left[\begin{array}{c} U_1 \\ U_2\end{array}\right]^{n+1}=J_2
  \left[\begin{array}{c} U_1 \\ U_2\end{array}\right]^n + \tilde F,
\end{equation}
where the associated iteration matrix is
\begin{equation}\label{eq:multi2sd3gen}
    J_2 =\left[\begin{array}{cc}
      (1+\sigma_{1})\Id-\mathbb{P}_1 & 0\\
      0 & (1+\sigma_2)\Id -\mathbb{P}_2
    \end{array}\right]^{-1} \left[\begin{array}{cc}
      0 & \sigma_{1} X \\
      \sigma_2 X &0
    \end{array}\right].
\end{equation}
Using the explicit formulas {\color{black}(\ref{eq:calp})} for the Calder\'on projectors, 
we can compute explicitly
\begin{equation}\label{eq:multi2sd3}
    J_2 = \displaystyle\left[\begin{array}{cccc}  0& 0&
     \frac{2\sigma_{1}+1}{2(\sigma_{1}+1)} & -  \frac{1}{2a(\sigma_{1}+1)} \\
    0& 0&  \frac{a}{2(\sigma_{1}+1)} & -\frac{2\sigma_{1}+1}{2(\sigma_{1}+1)} \\
    \frac{2\sigma_{2}+1}{2(\sigma_{2}+1)} & -  \frac{1}{2a(\sigma_{2}+1)} & 0&
     0\\
    \frac{a}{2(\sigma_{2}+1)} & -\frac{2\sigma_{2}+1}{2(\sigma_{2}+1)} & 0& 0 
    \end{array}\right],
\end{equation}
and the right hand side function is
\begin{equation}\label{Ftilde}
  \tilde{F}=\left[\begin{array}{c}
    -\frac{a\alpha(2\sigma_1+1)-\beta}{2a(1+\sigma_1)}\\
    -\frac{a\alpha-\beta(2\sigma_1+1)}{2(1+\sigma_1)}\\
    \frac{a\alpha(2\sigma_2+1)+\beta}{2a(1+\sigma_2)}\\
    \frac{a\alpha+\beta(2\sigma_2+1)}{2(1+\sigma_2)}
    \end{array}\right].
\end{equation}
The convergence factor of the block Jacobi iteration
\eqref{eq:blockJac} is given by the spectral radius of the iteration
matrix $J_2$, whose spectrum can be easily computed,
\begin{equation}
  \sigma(J_2)=\left \{-\sqrt{\frac{\sigma_{1}}{\sigma_{1}+1}},
    \sqrt{\frac{\sigma_{1}}{\sigma_{1}+1}}, -\sqrt{\frac{\sigma_{2}}{\sigma_{2}+1}},
    \sqrt{\frac{\sigma_{2}}{\sigma_{2}+1}} \right \}. 
\end{equation}
We note that the eigenvalues are independent of the problem parameter
$a$ and thus the convergence speed of the method only depends on the
relaxation parameters $\sigma_j$. This implies that the convergence
would be independent of the Fourier variable and thus robust when the
mesh size is refined in a two dimensional setting, as it was pointed
out in \cite{Dolean:2015:MFA}. Plotting the modulus of the eigenvalues
as function of $\sigma_j$, we obtain the result in Figure
\ref{fig:conv1d}.
We see that the algorithm diverges for $\sigma_j<-0.5$, stagnates for
$\sigma_j=-0.5$ and converges for all others values of $\sigma_j$.
For $\sigma_j$ close to zero, convergence is very rapid, and for
$\sigma_j=0$, $j=1,2$, the spectral radius of the iteration matrix
vanishes, which would make the method a direct solver, since the
iteration matrix becomes nil-potent. We have seen however also that for
vanishing relaxation parameters, the multitrace formulation
(\ref{eq:multi2sd1}) does not make sense any more, since the jump data
is not contained any more in the formulation. Nevertheless, the
associated block Jacobi iteration for the multitrace formulation
(\ref{eq:multi2sd1}) is well defined in the limit as $\sigma_j$ goes
to zero for $j=1,2$, and we get from \eqref{eq:multi2sd3}
\begin{equation}\label{eq:J2}
  \lim_{\sigma_j=0}J_2=\displaystyle\left[\begin{array}{cccc}  0& 0&
    \frac{1}{2} & -  \frac{1}{2a} \\
    0& 0&  \frac{a}{2} & -\frac{1}{2} \\
    \frac{1}{2} & -  \frac{1}{2a} & 0&
    0\\
   \frac{a}{2} & -\frac{1}{2} & 0& 0 
   \end{array}\right] =\left[\begin{array}{cc} 
   0 & \mathbb{P}X \\
   \mathbb{P}X & 0
  \end{array}\right].
\end{equation}
The limit of the right hand side (\ref{Ftilde}) is also well
defined, containing the data\footnote{Since we introduced only one
  multitrace $[\alpha,\beta]^T$ with jumps oriented from $\R^+$ to
  $\R^-$, in the first right hand side the operation
  $-X[\alpha,\beta]^T$ appears naturally to produce the consistent
  multitrace with the other orientation.}
$$
  \lim_{\sigma_j=0}\tilde{F}=\left[\begin{array}{c}
    -\frac{a\alpha-\beta}{2a}\\
    -\frac{a\alpha-\beta}{2}\\
    \frac{a\alpha+\beta}{2a}\\
    \frac{a\alpha+\beta}{2}
    \end{array}\right]
    =
    \left[\begin{array}{c}
      -\mathbb{P}X\left[\begin{array}{c}
        \alpha\\
        \beta
    \end{array}\right]\\
      \mathbb{P}\left[\begin{array}{c}
        \alpha\\
        \beta
    \end{array}\right]
  \end{array}\right].
$$
Therefore the block Jacobi iteration is also well defined in the limit
$\sigma_j=0$,
\begin{equation}\label{eq:opt}
  \left[\begin{array}{c} 
     U_1 \\ 
     U_2
  \end{array}\right]^{n+1}=\left[\begin{array}{cc} 
     0 & \mathbb{P}X \\
     \mathbb{P}X & 0
  \end{array}\right] \left[\begin{array}{c} 
    U_1 \\
    U_2
  \end{array}\right]^{n} + \left[\begin{array}{c} 
   -\mathbb{P}X[\alpha,\beta]^T  \\
   \mathbb{P}[\alpha,\beta]^T 
  \end{array}\right],
\end{equation}
and this iteration defines us at convergence a {\em new multitrace
  formulation}
\begin{equation}\label{eq:optlimit}
  \begin{array}{rcl} 
     U_1 -\mathbb{P}X U_2 &=& -\mathbb{P}X[\alpha,\beta]^T,\\ 
     U_2 -\mathbb{P}X U_1 &=& \mathbb{P}[\alpha,\beta]^T.
  \end{array}
\end{equation}
The advantage of this multitrace formulation is that it is already
preconditioned, block Jacobi applied to it is {\it optimal} in the
sense that convergence is achieved in a finite number of steps. A
direct calculation shows that $J_2^2$ equals zero, and thus
convergence is achieved in at most 2 iterations. We will see in
Section \ref{MTFOptSchwarzSec} that this iteration corresponds to a
well-known algorithm in domain decomposition.

\subsection{Multitrace Formulation for 3 Subdomains in 1D}
\label{ThreeDomainDDM}

We consider now a decomposition into three subdomains:
$I_1=(-\infty,-1)$, $I_0=(-1,1)$ and $I_2=(1,\infty)$, and functions
that satisfy
\begin{equation}
  \left\{\begin{array}{rcl}
    \displaystyle-\frac{d^2u}{dx^2} + a^2u&=&0,\, \mathbb{R}\setminus
    \{\pm 1\},\\
    \displaystyle\lim_{x\rightarrow\infty}|u(x)|&=&0.
  \end{array}\right.
\end{equation}
If we denote the restriction of the solution onto the subdomains by
$u_j = u|_{I_j}$, $j=0,1,2$, then by using a similar reasoning as in
the two-subdomain case in Subsection \ref{ReprForm1D}, we obtain the
representation formula
\begin{equation}\label{eq:repres3}
  \begin{array}{rcl}
u(x) &=& \displaystyle\left[-\frac{du_0}{dx}(-1)+ \frac{du_1}{dx}(-1)\right]{\cal
  G}(x+1) +\left[\frac{du_0}{dx}(1)- \frac{du_2}{dx}(1)\right]{\cal
  G}(x-1) \\
&-&\displaystyle [u_0(-1)-u_1(-1)]\frac{d{\cal
    G}}{dx}(x+1)+[u_0(1)-u_2(1)]\frac{d{\cal G}}{dx}(x-1)\\
&=& \displaystyle \left[-\frac{du}{dx}\right]_{-1}
\frac{e^{-a|x+1|}}{2a} +\left[\frac{du}{dx}\right]_{1}
\frac{e^{-a|x-1|}}{2a} \\
&+& \displaystyle [u]_{-1} \text{sign}(x+1)\frac{e^{-a|x+1|}}{2} -[u]_{1} \text{sign}(x-1)\frac{e^{-a|x-1|}}{2},
\end{array}
\end{equation}
where we defined the jumps of the derivatives to be
$$
  \beta_{-1}:=\left[-\frac{du}{dx}\right]_{-1} 
    := -\frac{du_0}{dx}(-1)+\frac{du_1}{dx}(-1),\quad 
   \beta_1:=\left[\frac{du}{dx}\right]_{1} 
    := \left[\frac{du_0}{dx}(1)- \frac{du_2}{dx}(1)\right],
$$
and the jumps in the function values are 
$$
  \alpha_{-1}:=[u]_{-1} :=  u_0(-1)-u_1(-1),\quad 
  \alpha_1:=[u]_{1}:= u_0(1)-u_2(1).
$$
Suppose now that we want to compute the Calder\'on projector for the
domain $I_0$. From \eqref{eq:repres3} we see that {\color{black}for $x\in I_{0}$ we have} 
\begin{equation}\label{eq:traces}
\begin{array}{rcl}
u_0(x) &=& \displaystyle \alpha_{-1} \frac{e^{-a(x+1)}}{2} + \beta_{-1}
\frac{e^{-a(x+1)}}{2a} + \alpha_{1}\frac{e^{a(x-1)}}{2}+\beta_{1}
\frac{e^{a(x-1)}}{2a} \\
& = & \displaystyle\displaystyle a \alpha_{-1} {\cal G}(x+1)+ \beta_{-1}
{\cal G}(x+1)+  a \alpha_{1}{\cal G}(x-1)+\beta_{1} {\cal G}(x-1), \\[2ex]
\displaystyle\frac{du_0}{dx}(x) & = &  \displaystyle -a\alpha_{-1} \frac{e^{-a(x+1)}}{2} -\beta_{-1}
\frac{e^{-a(x+1)}}{2} + a\alpha_{1}\frac{e^{a(x-1)}}{2}+\beta_{1}
\frac{e^{a(x-1)}}{2} \displaystyle \\
& =& \displaystyle\displaystyle -a^2 \alpha_{-1} {\cal G}(x+1)-a \beta_{-1}
{\cal G}(x+1)+  a^2 \alpha_{1}{\cal G}(x-1)+ a\beta_{1} {\cal G}(x-1). 
\end{array}
\end{equation}
If we define the Cauchy trace by
$$
  T_0(u) = \left[ u_0(-1),  -\frac{du_0}{dx}(-1), u_0(1), 
    \frac{du_0}{dx}(1) \right]^T,
$$
then from the formula \eqref{eq:traces} we obtain
$$
\begin{array}{rcl}
u_0(-1) & = & \displaystyle \alpha_{-1}\frac{1}{2} + \beta_{-1}
\frac{1}{2a} + a \alpha_{1}{\cal G}(-2)+\beta_{1} {\cal
  G}(-2),\\
\displaystyle -\frac{du_0}{dx}(-1) & =& \displaystyle \alpha_{-1}\frac{a}{2} + \beta_{-1}
\frac{1}{2} - a^2 \alpha_{1}{\cal G}(-2)-a\beta_{1} {\cal
  G}(-2),\\
u_0(1) & = & \displaystyle a\alpha_{-1}{\cal
  G}(2) + \beta_{-1}{\cal G}(2) +  \alpha_{1}\frac{1}{2}+\beta_{1} \frac{1}{2a},\\
\displaystyle \frac{du_0}{dx}(1) & =& \displaystyle-a^2\alpha_{-1}{\cal
  G}(2) - a\beta_{-1}{\cal G}(2) +  \alpha_{1}\frac{a}{2}+\beta_{1} \frac{1}{2},
\end{array}
$$
and thus using the short hand notation $g_\pm := {\cal G}(\pm 2)$
\begin{equation}
T_0(u) = \left[\begin{array}{cccc} \frac{1}{2} &
    \frac{1}{2a} & a g_-& g_- \\
\frac{a}{2} & \frac{1}{2} & -a^2 g_- & -ag_- \\
a g_+& g_+ & \frac{1}{2} &  \frac{1}{2a} \\
-a^2 g_+& -ag_+ & \frac{a}{2} &  \frac{1}{2} 
 \end{array}\right] 
\left(\begin{array}{c}  
\alpha_{-1} \\ \beta_{-1} \\  \alpha_{1} \\ \beta_{1}   
\end{array}  \right) =: \mathbb{P}_0 \left(\begin{array}{c}  
\alpha_{-1} \\ \beta_{-1} \\  \alpha_{1} \\ \beta_{1}   
\end{array}  \right).
\end{equation}
Here $\mathbb{P}_0$ is the Calder\'on projector for the middle subdomain,
\begin{equation} \label{eq:cal0}
\mathbb{P}_0 = \left [ \begin{array}{cc} \mathbb{P} & 2a\,g_- R \\ 2a\,g_+ R &
    \mathbb{P}\end{array}\right],\qquad R := \left[\begin{array}{cc} \frac{1}{2} &
    \frac{1}{2a} \\ -\frac{a}{2} & -\frac{1}{2}\end{array} \right],
\end{equation}
where $\mathbb{P}$ is given by the formula \eqref{eq:xp}. From the
facts that $\mathbb{P}$ is a projector, $\mathbb{P} R=0$,
$R\mathbb{P} =R$ and {\color{black}$R^{2} = 0$}, we see that $\mathbb{P}_0^2=\mathbb{P}_0$ and thus
$\mathbb{P}_0$ is a projector too, as expected. For the domains $I_1$
and $I_2$ with similar computations and definitions for the traces, we
obtain that $\mathbb{P}_1=\mathbb{P}_2=\mathbb{P}$.

The {\em multitrace formulation} in this case with three subdomains
states that the pairs $U_{1}$, $U_{2}$, and the quadruple
$U_0=(U_{01}^T,U_{02}^T)^T$ are traces of the solution defined on
$\Omega_j$ if they verify the relations
\begin{equation}\label{eq:multi3sd1}
\left\{\begin{array}{l}
(\Id-\mathbb{P}) U_1
+\sigma_1  \left(U_1-XU_{01}\right)=F_1,\\
{\color{black} (\Id-\mathbb{P}) U_{01}-2ag_-R U_{02}+
\sigma_0  \left(U_{01}-X U_1 \right) =F_{01},}\\
{\color{black} (\Id-\mathbb{P}) U_{02}-2ag_+R U_{01}+
\sigma_0  \left(U_{02}-X U_2 \right)=F_{02},}\\
(\Id-\mathbb{P}) U_2
+\sigma_2  \left(U_2-XU_{02}\right)= {\color{black} F_2 },
\end{array}\right.
\end{equation}
where $\sigma_{0,1,2}$ are again relaxation parameters. { \color{black} The right-hand side admits the 
explicit expression $F_1 = \sigma_1\lbr -\alpha_{-1},\beta_{-1}\rbr^{T}$,  
$F_2 = \sigma_2\lbr -\alpha_{1},\beta_{1}\rbr^{T}$ and 
$F_{01} = \sigma_0\lbr \alpha_{-1},\beta_{-1}\rbr^{T}$, 
$F_{02} = \sigma_0\lbr \alpha_{1},\beta_{1}\rbr^{T}$. }
In matrix form we obtain
\begin{equation}\label{eq:multi3sd2}
\left[\begin{array}{cccc}
(1+\sigma_1)\Id-\mathbb{P} & -\sigma_1 X & 0 & 0\\
-\sigma_0 X & (1+\sigma_0)\Id-\mathbb{P} & -2a\,g_-R & 0 \\
0 & -2a\,g_+R &  (1+\sigma_0)\Id-\mathbb{P} & -\sigma_0 X \\
0& 0& -\sigma_2 X & (1+\sigma_2)\Id-\mathbb{P} 
\end{array}\right] \left[\begin{array}{c} U_1 \\ {\color{black}U_{01}}\\ {\color{black}U_{02}} \\
  U_2\end{array}\right] = \left[\begin{array}{c} F_1 \\ {\color{black}F_{01}}\\ {\color{black}F_{02}}\\ {\color{black} F_2 }
  \end{array}\right].
\end{equation}
As in the case of two subdomains, it is natural to apply a
block-Jacobi iteration to \eqref{eq:multi2sd2}, which leads to the
iteration
\begin{equation}\label{eq:block}
\left[\begin{array}{c} U_1 \\ U_0 \\
  U_2\end{array}\right]^{n+1} = J_3 \left[\begin{array}{c} U_1 \\ U_0 \\
  U_2\end{array}\right]^{n} + \tilde F, 
 \quad \tilde F = 
\left[\begin{array}{l}
(1+\sigma_{1})^{-1}(\sigma_{1}\Id+\mathbb{P})\lbr -\alpha_{-1},\beta_{-1}\rbr^{T}\\
(1+\sigma_{0})^{-1}(\sigma_{0}\Id+\mathbb{P})\lbr \phantom{-}\alpha_{-1},\beta_{-1}\rbr^{T}\\
(1+\sigma_{0})^{-1}(\sigma_{0}\Id+\mathbb{P})\lbr \phantom{-}\alpha_{1},\beta_{1}\rbr^{T}\\
(1+\sigma_{2})^{-1}(\sigma_{2}\Id+\mathbb{P})\lbr -\alpha_{1},\beta_{1}\rbr^{T}
\end{array}\right], 
\end{equation}
where the iteration matrix is given by
\begin{equation}\label{eq:multi3sd3}
\begin{array}{rcl}
J_3 &=&\left[\begin{array}{cccc}
(1+\sigma_1)\Id-\mathbb{P} & 0 & 0 & 0\\
0 & (1+\sigma_0)\Id -\mathbb{P} & 0 & 0 \\
0 & 0& (1+\sigma_0)\Id -\mathbb{P} & 0  \\
0 & 0& 0& (1+\sigma_2)\Id -\mathbb{P} 
\end{array}\right]^{-1} \\[3ex]
& \cdot & \left[\begin{array}{cccc}
0& \sigma_1 X & 0 & 0\\
\sigma_0 X & 0 & 2a\,g_-R & 0 \\
0 & 2a\,g_+R &  0 & \sigma_0 X \\
0& 0& \sigma_2 X & 0
\end{array}\right]. 
\end{array}
\end{equation}
The convergence factor of the block Jacobi iteration is again
determined by the eigenvalues of the iteration matrix $J_3$, which are
readily calculated to be
\begin{equation}
\sigma(J_3)=\left \{-\sqrt{\frac{\sigma_j}{\sigma_j+1}},
    \sqrt{\frac{\sigma_j}{\sigma_j+1}}, j=0,1,2\right \}. 
\end{equation}
We see that the convergence behavior with three subdomains is
identical to the case of two subdomains, and in the limiting case
when $\sigma_j=0$, we obtain for the limit of the iteration
$J_3$
\begin{equation}\label{eq:multi3}
J_3 =\left[\begin{array}{cccc}
0 & \mathbb{P}X & 0& 0\\
\mathbb{P}X & 0& 2ag_- R & 0 \\
0 & 2ag_+ R &0 & \mathbb{P}X \\
0 & 0&\mathbb{P}X & 0
\end{array}\right].
\end{equation}
In this case it is easy to check that $J_3^4=0$, and therefore algorithm
\eqref{eq:block} converges in at most $4$ iterations.

\section{Multitrace Formulations and Optimal Schwarz Methods}\label{MTFOptSchwarzSec}

We now want to relate the block Jacobi iteration we defined for the
multitrace formulation (\ref{eq:multi2sd1}) to a well studied class of
domain decomposition methods of Schwarz type. {\color{black} While the
  analysis of this section also holds for Problem (\ref{ModelEq2}),
  the form of the associated Calder\'on projectors (\ref{eq:calp}) has
  become too simple because of the strong symmetries to find the
  relation between the multitrace formulation and optimal Schwarz
  methods. We thus first have to study the Calder\'on projectors for a
  non-symmetric domain configuration on a bounded domain.}

%and we thus study first the same case on a bounded
%domain.

\subsection{Calder\'on Projectors on a Bounded Domain}

We are interested in the solution of the transmission problem
\begin{equation}\label{transmissionproblem}
\left\{\begin{array}{rcl}
-u''(x)+a^2u(x) &=& 0, x\in (0,1)\setminus \{\gamma\},\\
\displaystyle [u] &=& \alpha,\\
\displaystyle [u'] & = & \beta,\\
u(0)&=&u(1)=0.
\end{array}\right.
\end{equation}
Local solutions to the left and right of the jumps satisfying the
outer boundary conditions are given by
\begin{equation}
  u_1(x):=u|_{(0,\gamma)} = c_1 \sinh(ax),\quad
    u_2(x):=u|_{(\gamma,1)} = c_2 \sinh(a(1-x)),
\end{equation}
where $c_j$, $j=1,2$ are some constants.  Using the same expressions
for the jumps as in the unbounded case,
$$
u_2(\gamma)-u_1(\gamma)=\alpha,\quad -u_2'(\gamma)+u_1'(\gamma)=\beta,
$$
we obtain an equation for the constants $c_1$ and $c_2$,
$$
\left[\begin{array}{cc}
\sinh(a(1-\gamma)) & -\sinh(a\gamma) \\
a\cosh(a(1-\gamma)) & a \cosh(a\gamma)
\end{array}\right]\left[\begin{array}{c} c_2 \\ c_1\end{array}\right] = \left[\begin{array}{c} \alpha \\ \beta\end{array}\right].
$$
Solving the linear system for the constants leads then to the closed
form solutions of the transmission problem (\ref{transmissionproblem}),
\begin{equation}
\left\{
\begin{array}{rcl}
u_1(x) & = & \displaystyle\frac{1}{D} \left[-a\cosh(a(1-\gamma))\alpha
  +\sinh(a(1-\gamma))\beta\right]\sinh(ax), \\[1ex]
u_2(x) & = & \displaystyle\frac{1}{D}
\left[a\cosh(a\gamma)\alpha+\sinh({\color{black}a}\gamma)\beta \right]\sinh(a(1-x)),
\end{array}
\right.
\end{equation}
where $
D:=a\left[\cosh(a(1-\gamma))\sinh({\color{black}a}\gamma)+\sinh(a(1-\gamma))
  \cosh({\color{black}a}\gamma)\right]$. Proceeding as in the unbounded case, we
can deduce that if $u_2$ is a function satisfying the equation on
$(\gamma,1)$, then it can be expressed as $u_2(x) = (G_{2}\circ
T_{2}(u))(x)$, where
\begin{equation}\label{eq:trace2}
  \arraycolsep0.1em
  G_2\left(\begin{array}{c}\alpha \\
    \beta\end{array}\right):=\frac{1}{D}\left[\alpha\cosh(a\gamma)
      +\beta\sinh(a\gamma)\right]\sinh(a(1-x)),\ 
  T_2(u):=\left(\begin{array}{c} u(\gamma_+) \\
    -\frac{du}{dx}(\gamma_+) \end{array}\right). 
\end{equation}
Again $T_{2}\circ G_{2}$ is a $2\times 2$ matrix whose coefficients can be 
explicitly computed, 
\begin{equation}\label{eq:p2}
  \mathbb{P}_2:= T_2\circ G_2 = \frac{1}{D} \left[\begin{array}{cc}
    a\cosh(a\gamma)\sinh(a(1-\gamma)) & \sinh(a\gamma)\sinh(a(1-\gamma)) \\
    a^2\cosh(a\gamma)\cosh(a(1-\gamma)) & {\color{black}a}\sinh(a\gamma)\cosh(a(1-\gamma)) 
  \end{array}\right].
\end{equation}
With a similar reasoning on $(0,\gamma)$ we obtain
\begin{equation}\label{eq:p1}
\mathbb{P}_1:= T_1\circ G_1 = \frac{1}{D} \left[\begin{array}{cc}
a\cosh(a(1-\gamma))\sinh(a\gamma) & \sinh(a(1-\gamma))\sinh(a\gamma) \\
a^2\cosh(a(1-\gamma))\cosh(a\gamma) & {\color{black}a}\sinh(a(1-\gamma))\cosh(a\gamma) 
\end{array}\right],
\end{equation}
where 
\begin{equation}\label{eq:trace1}
  \arraycolsep0.1em
  G_1\left(\begin{array}{c}\alpha \\
    \beta\end{array}\right):=\frac{1}{D}\left[\alpha\cosh(a(1-\gamma))
      +\beta\sinh(a(1-\gamma))\right]\sinh(ax),\ 
  T_1(u):=\left(\begin{array}{c} u(\gamma_-) \\
    \frac{du}{dx}(\gamma_-) \end{array}\right). 
\end{equation}
As in the unbounded domain case, the two operators $\mathbb{P}_{1,2}$
are projectors, $\mathbb{P}_{j}^2=\mathbb{P}_{j}$, and they are called
{\em Calder\'on projectors}.

If we consider the same multitrace formulation (\ref{eq:multi2sd1}) as
in the unbounded case and apply a block-Jacobi iteration, we obtain
for the iteration matrix in an analogous way
\begin{equation}\label{eq:multi2sdgen}
  \begin{array}{rcl}
    J_2 &=&\left[\begin{array}{cc}
     (1+\sigma_{1})\Id-\mathbb{P}_1 & 0\\
    0 & (1+\sigma_2)\Id -\mathbb{P}_2
    \end{array}\right]^{-1} \left[\begin{array}{cc}
    0 & \sigma_{1} X \\
   \sigma_2 X &0
\end{array}\right]\\[2ex]
& =& \displaystyle\left[\begin{array}{cc}
0 & \displaystyle\frac{1}{\sigma_1+1}(\sigma_{1} I+\mathbb{P}_{1})X \\
\displaystyle\frac{1}{\sigma_1+1}(\sigma_{2} I+\mathbb{P}_{2})X  &0
\end{array}\right],
\end{array}
\end{equation} 
where the second equality holds since the $\mathbb{P}_{j}$ are
projectors, and we hence do not need to rely on explicit expressions
to obtain this result! We thus obtain an identical convergence
behavior like in the unbounded domain case and in the limiting case
$\sigma_j=0$ the optimal iteration
\begin{equation}\label{eq:optgen}
\left[\begin{array}{c} U_1 \\ U_2\end{array}\right]^{n+1}= \left[\begin{array}{cc} 
0 & \mathbb{P}_{1}X \\
\mathbb{P}_{2}X & 0
\end{array}\right] \left[\begin{array}{c} U_1 \\
  U_2\end{array}\right]^{n} + {\color{black}
 \left[\begin{array}{c}
      -\mathbb{P}_1X\left[\begin{array}{c}
        \alpha\\
        \beta
    \end{array}\right]\\
      \mathbb{P}_2\left[\begin{array}{c}
        \alpha\\
        \beta
    \end{array}\right]
  \end{array}\right]},
\end{equation}
with a right hand side corresponding to the bounded domain case.

\subsection{Dirichlet to Neumann Operators and Calder\'on Projectors}
\label{DtNAndCalderonSec}

To find a relation between the optimal block Jacobi iteration for the
multitrace formulation and the optimal Schwarz methods, we now write
the Calder\'on projectors in terms of the Dirichlet to Neumann (DtN)
operators. First we compute the DtN operators on the domains
$\Omega_1=(0,\gamma)$ and $\Omega_2= (\gamma,1)$. To start, we consider
the boundary value problem
\begin{equation}
\left\{\begin{array}{rcl}
-u_1''(x)+a^2u_1(x) &=& 0, x\in (0,\gamma),\\
u_1(\gamma)& = & g,\\
u_1(0)&=&0.
\end{array}\right.
\end{equation}
Then the $DtN_1$ associates to the Dirichlet data $g=u_1(\gamma)$ the
normal derivative of the solution $u_1'(\gamma)$. A simple computation
gives
$$
  u_1(x) = \frac{\sinh(ax)}{\sinh(a\gamma)}g\quad \Longrightarrow\quad 
    u_1'(\gamma) = \frac{a \cosh(a\gamma)}{\sinh(a\gamma)}g =:\mbox{DtN}_1 g.
$$
We consider next the boundary value problem 
\begin{equation}
\left\{\begin{array}{rcl}
-u_2''(x)+a^2u_2(x) &=& 0, x\in (\gamma,1),\\
u_2(\gamma)& = & g,\\
u_2(1)&=&0.
\end{array}\right.
\end{equation}
Then the $DtN_2$ associates to the Dirichlet data $g=u_2(\gamma)$ the
normal derivative of the solution $-u_2'(\gamma)$, and we obtain by a
direct calculation
$$
  u_2(x) = \frac{\sinh(a(1-x))}{\sinh(a(1-\gamma))}g\quad
  \Longrightarrow \quad -u_2'(\gamma) =
    \frac{a \cosh(a(1-\gamma))}{\sinh(a(1-\gamma))}g =:\mbox{DtN}_2 g.
$$
Similarly we can define the Neumann to Dirichlet operators
$\mbox{NtD}_j$, which calculate from given Neumann data the associated
Dirichlet data, and are thus just the inverses of the corresponding
$\mbox{DtN}_j$. 

Comparing the expressions for the $\mbox{DtN}_j$ and $\mbox{NtD}_j$
operators with the expressions for the Calder\'on projectors in
(\ref{eq:p2}) and (\ref{eq:p1}), we see that the Calder\'on projectors
can be re-written as
\begin{equation}\label{eq:p12}
\begin{array}{rcl}
\mathbb{P}_1 &=&
\left[\begin{array}{cc} (\mbox{DtN}_1+\mbox{DtN}_2)^{-1}\mbox{DtN}_2 &  (\mbox{DtN}_1+\mbox{DtN}_2)^{-1}\\
(\mbox{NtD}_1+\mbox{NtD}_2)^{-1}& (\mbox{NtD}_1+\mbox{NtD}_2)^{-1}\mbox{NtD}_2 
\end{array}\right],\\[2ex]
\mathbb{P}_2 &=&
\left[\begin{array}{cc} (\mbox{DtN}_1+\mbox{DtN}_2)^{-1}\mbox{DtN}_1 &  (\mbox{DtN}_1+\mbox{DtN}_2)^{-1}\\
(\mbox{NtD}_1+\mbox{NtD}_2)^{-1}& (\mbox{NtD}_1+\mbox{NtD}_2)^{-1}\mbox{NtD}_1 
\end{array}\right].
\end{array}
\end{equation}
This reformulation of the Calder\'on operators allows us in the next
section to identify the optimal block Jacobi method with a well
understood optimal Schwarz method.

\subsection{Relation to Optimal Schwarz Methods}

Let ${\cal L}:=-\partial_{xx}+a^2$ be the differential operator we
have been studying so far. A non-overlapping optimal Schwarz iteration
(see \cite{gander2006optimized} and references therein)
using the decomposition into the two subdomains $\Omega_1=(0,\gamma)$
and {\color{black}$\Omega_2=(\gamma,1)$} from Subsection \ref{DtNAndCalderonSec} is
given by the algorithm
\begin{equation}\label{OptimalSchwarz}
\begin{array}{rcll}
{\cal L}u_1^{n+1}&=&f & \mbox{in $\Omega_1$}\\
\displaystyle\frac{\partial u_1^{n+1}}{\partial x} +\mbox{DtN}_2
u_1^{n+1} &=& \displaystyle\frac{\partial u_2^{n}}{\partial x}
+\mbox{DtN}_2 u_2^{n} & \mbox{on $x=\gamma$},\\
{\cal L}u_2^{n+1}&=&f & \mbox{in $\Omega_1$},\\
\displaystyle-\frac{\partial u_2^{n+1}}{\partial x} +\mbox{DtN}_1
u_2^{n+1} &=&\displaystyle-\frac{\partial u_1^{n}}{\partial x}
+\mbox{DtN}_1 u_1^{n}& \mbox{on $x=\gamma$}.
\end{array}
\end{equation}
It is well known, see for example \cite{gander2006optimized}, that the
optimal Schwarz algorithm (\ref{OptimalSchwarz}) converges in two
iterations, like the block-Jacobi algorithm \eqref{eq:blockJac} with
two subdomains and relaxation parameter $\sigma_j=0$, $j=1,2$.
Schwarz methods are however usually not used to solve transmission
problems, and zero jumps are enforced by the algorithm
\eqref{OptimalSchwarz} at the interface $\gamma$. To study the
convergence of algorithm \eqref{OptimalSchwarz}, one analyzes directly the
error equations, i.e. algorithm (\ref{OptimalSchwarz}) with right hand
side $f=0$, and studies how the subdomain iterates go to zero as the
iteration progresses. In this homogeneous case, the iterates
$u_j^{n+1}$, $j=1,2$ are solutions of the homogeneous problems inside
the subdomains, and the normal derivatives can be expressed in terms
of the $\mbox{DtN}$ operators: for example $\frac{\partial
  u_1^{n+1}}{\partial x}=\mbox{DtN}_1u_1^{n+1}$ on $x=\gamma$. This
means that the iteration on the the first subdomain can be written
directly on the interface $x=\gamma$ as a function of the Dirichlet
trace of the iterate,
\begin{equation}\label{DirTrace}
\begin{array}{rcl}
(\mbox{DtN}_1+\mbox{DtN}_2)u_1^{n+1}&=&\displaystyle\frac{\partial
  u_2^{n}}{\partial x} +\mbox{DtN}_2 u_2^{n}, 
\\[2ex]
\Longleftrightarrow \quad u_1^{n+1} & = & \displaystyle (\mbox{DtN}_1+\mbox{DtN}_2)^{-1}\left(\frac{\partial
  u_2^{n}}{\partial x} +\mbox{DtN}_2 u_2^{n} \right).
\end{array}
\end{equation}
It is also possible to write this iteration based on the Neumann
traces, namely
\begin{equation}\label{NeuTrace}
\begin{array}{rcl}
\displaystyle\frac{\partial u_1^{n+1}}{\partial x} +\mbox{DtN}_2 \mbox{NtD}_1
\frac{\partial u_1^{n+1}}{\partial x} &=&\displaystyle\frac{\partial
  u_2^{n}}{\partial x} +\mbox{DtN}_2 u_2^{n}, \\
\Longleftrightarrow\quad \displaystyle\frac{\partial u_1^{n+1}}{\partial x} & = & \displaystyle (\mbox{DtN}_2 \mbox{NtD}_1+I)^{-1}\left(\frac{\partial
  u_2^{n}}{\partial x} +\mbox{DtN}_2 u_2^{n} \right), 
\\[2ex]
\Longleftrightarrow\quad \displaystyle\frac{\partial u_1^{n+1}}{\partial x} & = & \displaystyle (\mbox{NtD}_1+\mbox{NtD}_2)^{-1}\left(\mbox{NtD}_2\frac{\partial
  u_2^{n}}{\partial x} + u_2^{n} \right),
\end{array}
\end{equation}
where we used that $\mbox{DtN}_j$ is the inverse of the
$\mbox{NtD}_j$. Combining the two formulations (\ref{DirTrace}) and
(\ref{NeuTrace})\footnote{which means we would run the optimal Schwarz
  algorithm twice simultaneously, once on the Dirichlet traces and
  once on the Neumann traces, which would be very costly and not
  advisable in practice}, we obtain the iteration
\begin{equation}\label{OptSchwarzAsCalderon}
\left[\begin{array}{c}
u_1^{n+1}\\
\displaystyle\frac{\partial u_1^{n+1}}{\partial x}
\end{array}\right] = \left[\begin{array}{cc} (\mbox{DtN}_1+\mbox{DtN}_2)^{-1}\mbox{DtN}_2 &  (\mbox{DtN}_1+\mbox{DtN}_2)^{-1}\\
(\mbox{NtD}_1+\mbox{NtD}_2)^{-1}& (\mbox{NtD}_1+\mbox{NtD}_2)^{-1}\mbox{NtD}_2 
\end{array}\right] \left[\begin{array}{c}
u_2^{n}\\
\displaystyle\frac{\partial u_2^{n}}{\partial x}
\end{array}\right], 
\end{equation}
and we see the first Calder\'on projector $\mathbb{P}_1$ appear from
(\ref{eq:p12}). By re-writing this relation in terms of the traces
from Subsection \ref{2Dom1D} and taking into account the sign
convention we used there, iteration (\ref{OptSchwarzAsCalderon}) is
identical to
\begin{equation}\label{OptSchwarzCalderonHom}
  U_1^{n+1} = \mathbb{P}_1X U_2^n,\quad \mbox{and similarly}\quad
  U_2^{n+1} = \mathbb{P}_2X U_1^n,
\end{equation}
which is obtained similarly for the second subdomain. By comparing
with \eqref{eq:optgen}, we see that iteration
\eqref{OptSchwarzCalderonHom} is identical to \eqref{eq:optgen} in the
homogeneous case, i.e. when the jumps are zero. We have thus proved
the following
\begin{theorem}
  For two subdomains, the optimal multitrace iteration
  \eqref{eq:optgen} is an equivalent algorithm to the optimal Schwarz
  iteration \eqref{OptimalSchwarz}: it runs the optimal Schwarz
  algorithm twice simultaneously, once on the Dirichlet traces and
  once on the Neumann traces.
\end{theorem}

\section{General Multitrace Formulation}\label{MTFGeneral}

We now illustrate what insight can be gained from our simple problem
for multitrace formulations in a higher dimensional, geometrically more
general context using the common multitrace formalism. Although we do
not wish to dwell on the functional analytic aspects of boundary
integral equations, we need to introduce functional spaces adapted to
integral operators. Given a Lipschitz domain $\Omega\subset \R^{d}$,
we will consider the space of square integrable functions
$\mL^{2}(\Omega) = \{v\;, \;\Vert v\Vert_{\mL^{2}(\Omega)}^{2} =
\int_{\Omega} \vert v(\bfx)\vert d\bfx <+\infty\}$, and the Sobolev
spaces $\mH^{1}(\Omega):=\{ v\in \mL^{2}(\Omega), \;\nabla v\in
\mL^{2}(\Omega)\}$ and $\mH^{1}(\Delta, \Omega):=\{ v\in
\mH^{1}(\Omega),\;\Delta v\in \mL^{2}(\Omega)\}$ equipped with the
associated natural norms $\Vert v\Vert_{\mH^{1}(\Omega)}^{2} = \Vert
v\Vert_{\mL^{2}(\Omega)}^{2} + \Vert \nabla
v\Vert_{\mL^{2}(\Omega)}^{2}$, and $\Vert
v\Vert_{\mH^{1}(\Delta,\Omega)}^{2} = \Vert
v\Vert_{\mH^{1}(\Omega)}^{2} + \Vert \Delta
v\Vert_{\mL^{2}(\Omega)}^{2}$.

We also need to introduce trace spaces. First of all, the application
$v\mapsto v\vert_{\partial\Omega}$ extends to a continuous operator
mapping $\mH^{1}(\Delta, \Omega)$ to a strict subspace of
$\mL^{2}(\partial\Omega)$ that we denote by
$\mH^{1/2}(\partial\Omega):=\{ v\vert_{\partial\Omega},\; v\in
\mH^{1}(\Delta,\Omega)\;\}$, equipped with the norm $\Vert
v\Vert_{\mH^{1/2}(\partial\Omega)} = \inf\{\Vert
u\Vert_{\mH^{1}(\Omega)},\;u\vert_{\partial\Omega} = v\}$. Finally,
$\mH^{-1/2}(\partial\Omega)$ denotes the dual space to
$\mH^{-1/2}(\partial\Omega)$, equipped with the associated canonical
dual norm. Denoting by $\bfn$ the normal vector to $\partial\Omega$
pointing outward, it is a consequence of Rademacher's theorem that
the application $v\mapsto \bfn\cdot\nabla v\vert_{\partial\Omega}$ can
be extended as a continuous map of $\mH^{1}(\Delta,\Omega)$ onto
$\mH^{-1/2}(\partial\Omega)$, see \cite[Thm.2.7.7]{MR2743235}.

\subsection{Representation Formulas} 

We show now how the concrete representation formulas from the one
dimensional example of Subsection \ref{ReprForm1D} look for domains
$\Omega \subset \R^{d}$ with $d = 1,2,3,\dots$. Given $a>0$, we are
still interested in problems of the form $-\Delta u + a^2u=0$ in
domains of $\R^{d}$.  In what follows, $\mathcal{G}$ refers to the
unique Green kernel of this equation that decreases at infinity, i.e.
$$
  -\Delta \mathcal{G} + a^2\mathcal{G}=\delta_{0}(\bfx)\quad
    \textrm{in}\;\R^{d}\setminus\{0\},
    \quad\quad\lim_{\vert\bfx\vert\to \infty}\mathcal{G}(\bfx) = 0,
$$ 
where $\delta_{0}$ is the Dirac distribution centered at $\bfx =
0$. Explicit expressions of $\mathcal{G}$ (depending on the dimension
of the space) are known.  For $d=3$ for example, $\mathcal{G}(\bfx) =
\exp(-a\vert\bfx\vert)/(4\pi\vert\bfx\vert)$.  In this paragraph,
$\Omega\subset \R^{d}$ will refer to a Lipschitz open set with bounded
boundary, i.e. $\Omega$ is bounded or the complementary of a bounded
set. Associated to this domain, we consider the potential operator
\begin{equation}\label{DefPot}
G(v,q)(\bfx) := \int_{\partial\Omega}q(\bfy)\mathcal{G}(\bfx-\bfy)+v(\bfy)\,
\bfn(\bfy)\cdot(\nabla\mathcal{G})(\bfx-\bfy)\,d\sigma(\bfy).
\end{equation}
In this definition $\bfn$ refers to the normal vector field to 
$\partial\Omega$ pointing toward the exterior of $\Omega$. 
The potential operator $G$ maps continuously arbitrary pairs of traces 
$(v,q)\in\mH^{1/2}(\partial\Omega)\times \mH^{-1/2}(\partial\Omega)$
to functions $u = G(v,q)$ satisfying $-\Delta u +a^{2} u = 0$ in 
$\R^{d}\setminus\partial\Omega$. Analogous to (\ref{eq:tracep}), consider 
the trace operator 
\begin{equation}\label{DefTraceOp}
T(u) :=  \left( \begin{array}{r}
u\vert_{\partial\Omega}\\
\bfn\cdot \nabla u\vert_{\partial\Omega}
\end{array}\right).
\end{equation}
This definition makes sense for $u\in
\mH^{1}(\Delta,\Omega) = \{ u\in \mH^{1}(\Omega)\;\vert\;
\Delta u\in \mL^{2}(\Omega) \}$. We underline also that, in the
definition of $T$, the traces are taken from the \textit{interior} of
$\Omega$.  The next result is proved for example in \cite[Theorem
  3.1.6]{MR2743235}.\\
\begin{proposition}\label{PropReprFormula}
Let $u\in \mH^{1}(\Omega)$ satisfy $-\Delta u+a^{2}u = 0$ 
in $\Omega$. We have the representation formula 
\begin{equation}\label{ReprFormula}
G(T(u))(\bfx)\;=\;
\left\{\begin{array}{ll}
u(\bfx) & \textrm{for}\;\;\bfx\in \Omega,\\
0       & \textrm{for}\;\;\bfx\in \R^{d}\setminus\overline{\Omega}\;.
\end{array}\right.
\end{equation}
\end{proposition}

\subsection{Calder\'on Projectors}

For any pair of traces $V = (v,q)$, the function $u(\bfx) =
G(V)(\bfx)$ satisfies $-\Delta u + a^{2}u = 0$ in $\Omega$, so we can
apply the representation formula (\ref{ReprFormula}) above, like we
applied the representation formula in the one dimensional case in
Subsection \ref{Calderon1d}, which yields $G(T\cdot G(V))(\bfx) =
G(V)(\bfx)$ for $\bfx\in \Omega$. Taking the traces of this identity
leads to $(T\cdot G) (T\cdot G)(V) = (T\cdot G)(V)$. Setting
$\mathbb{P} := T\cdot G$, we see that $\mathbb{P}^{2} = \mathbb{P}$,
and hence the operator $\mathbb{P}$ is a projector, called {\em
  Calder\'on projector} associated to $\Omega$.

\subsection{Multitrace Formulation with 2 Subdomains}\label{2Dom2D}

We consider now a higher-dimensional counterpart of the one dimensional
two subdomain situation studied in Subsection \ref{2Dom1D}. Let
$\Omega_{1}\subset \R^{d}$ refer to any bounded Lipschitz subdomain
and set $\Omega_{2} := \R^{d}\setminus \overline{\Omega}_{1}$, $\Gamma
:= \partial\Omega_{1}$. In what follows we denote by $G_{j}$, $j=1,2$
the potential operator given by Formula (\ref{DefPot}) with $\Omega =
\Omega_{j}$ and $\bfn = \bfn_{j}$. The Calder\'on projector associated
to $\Omega_{j}$ will be denoted $\mathbb{P}_{j}$.

We first point out some remarkable identities relating
$\mathbb{P}_{1}$ to $\mathbb{P}_{2}$. First observe that, since
$\bfn_{2} = -\bfn_{1}$, we have $G_{2}(U) = -G_{1}( X U)$ for all
$U\in \mH^{1/2}(\Gamma)\times\mH^{-1/2}(\Gamma) $ where $X$ is the
matrix defined in (\ref{eq:xp}). Assume that $U = (\alpha,\beta)$, with
$\alpha\in\mH^{1/2}(\Gamma)$, $\beta\in \mH^{-1/2}(\Gamma)$, and consider
the unique function $u\in\mH^{1}(\R^{d}\setminus\Gamma)$ satisfying
$-\Delta u + a^{2}u = 0$ in $\R^{d}\setminus\Gamma$, $\lbr
u\rbr_{\Gamma}= (u\vert_{\Omega_{1}})\vert_{\Gamma} -
(u\vert_{\Omega_{2}})\vert_{\Gamma} = \alpha$, $\lbr
\partial_{n}u\rbr_{\Gamma} = \beta$, so that $U = T_{1}(u) -
XT_{2}(u)$. Applying Proposition \ref{PropReprFormula} both to
$u\vert_{\Omega_{1}}$ and $u\vert_{\Omega_{2}}$ yields
$$
\begin{array}{rl}
(T_{1} - XT_{2})G_{1}(U) 
& = (T_{1} - XT_{2})(\;G_{1}(T_{1}(u))-G_{1}(XT_{2}(u))\;)\\[5pt]
& = (T_{1} - XT_{2})(\;G_{1}(T_{1}(u))+G_{2}(T_{2}(u))\;)\\[5pt]
& = T_{1}\cdot G_{1}(T_{1}(u))  - XT_{2}\cdot G_{2}(T_{2}(u))\\[5pt] 
& = T_{1}(u) - XT_{2}(u) = U.

\end{array}
$$ 
Since $U$ was chosen arbitrarily, we conclude from this that $(T_{1} -
XT_{2})G_{1} = \Id$, and since $G_{1} = -G_{2}X$ we finally obtain the
identity
\begin{equation}\label{CommutingRelation}
  X\mathbb{P}_{2}X = \Id - \mathbb{P}_{1}.
\end{equation}
Now we want to consider a transmission problem similar to
(\ref{ModelEq2}).  Given boundary data $h
=(h_{\mathrm{D}},h_{\mathrm{N}}) \in \mH^{1/2}(\Gamma)\times
\mH^{-1/2}(\Gamma)$, we consider the transmission problem
\begin{equation}\label{ModelEq3}
\left\{\begin{array}{l}
u\in \mH^{1}(\R^{d}),\\
-\Delta u + a^{2} u = 0\quad\textrm{in}\;\; \R^{d}\setminus\Gamma,\\
\lbr u\rbr_{\Gamma} = h_{\mathrm{D}},\quad 
\lbr \partial_{n}u\rbr_{\Gamma} = h_{\mathrm{N}},
\end{array}\right.
\end{equation}
where $\lbr u\rbr_{\Gamma} := (u\vert_{\Omega_{1}})\vert_{\Gamma} -
(u\vert_{\Omega_{2}})\vert_{\Gamma}$ and $\lbr
\partial_{n}u\rbr_{\Gamma} := \bfn_{1}\cdot\nabla
(u\vert_{\Omega_{1}})\vert_{\Gamma} +
\bfn_{2}\cdot\nabla(u\vert_{\Omega_{2}})\vert_{\Gamma}$ for the
Dirichlet and Neumann jumps of the traces.  Setting $U_{1} :=
T_{1}(u)$ and $U_{2} := T_{2}(u)$, the jump conditions in the
equations above can be rewritten as $T_{1}(u) - X T_{2}(u) = h$. The
local multitrace formulation associated to Problem (\ref{ModelEq2})
is precisely of the same form as in the simple one dimensional case
(\ref{eq:multi2sd1}), namely
\begin{equation}\label{MultiTr2D1}
\left\{\begin{array}{l}
(\Id-\mathbb{P}_1) U_1
+\sigma_1  \left(U_1-X U_2 \right)=F_{1},\\[5pt]
(\Id-\mathbb{P}_2) U_2
+\sigma_2  \left(U_2-XU_1\right)=F_{2},
\end{array}\right.
\end{equation}
%\marginpar{Need to check $F_j$ here, since in 1d there was an error}
with {\color{black} $F_{1} = \sigma_{1}h$ and $F_{2} = -\sigma_{2}X h$}. This time
however, the operator associated to (\ref{MultiTr2D1}) is not a simple
$4\times 4$ matrix any more with complex scalar entries, it is a
$4\times 4$ matrix of integral operators. The block Jacobi iteration
operator associated with this formulation is
$$
\mathbb{J}_{2} =\left[\begin{array}{cc}
(1+\sigma_{1})\Id-\mathbb{P}_1 & 0\\
0 & (1+\sigma_2)\Id -\mathbb{P}_2
\end{array}\right]^{-1} \left[\begin{array}{cc}
0 & \sigma_{1} X \\
\sigma_2 X &0
\end{array}\right].
$$
To simplify this expression, note that for any $\gamma\in\C$ and
$j=1,2$, since $\mathbb{P}^{2}_{j} = \mathbb{P}_{j}$, we have the
identity
\begin{equation}\label{ExprInv}
  ((1+\gamma)\Id -\mathbb{P}_{j})(\gamma\Id +\mathbb{P}_{j}) 
    = \gamma(1+\gamma)\Id.
\end{equation}
Taking this identity into account with $\gamma = \sigma_{j}$, $j=1,2$ leads 
to a simplified expression of the Jacobi iteration matrix as in the
one dimensional case where we first used direct manipulations,
\begin{equation}\label{ExprJacMat2}
  \mathbb{J}_{2} =\left[\begin{array}{cc}
    0 & (1+\sigma_{1})^{-1}(\sigma_{1}\Id +\mathbb{P}_{1})X \\
    (1+\sigma_{2})^{-1}(\sigma_{2}\Id +\mathbb{P}_{2}) X & 0
  \end{array}\right].
\end{equation}
To compute the eigenvalues of this operator, it is convenient to
first square it. As a preliminary remark note that, according to
(\ref{CommutingRelation}) and since $\mathbb{P}_{1}^{2} =
\mathbb{P}_{1}$, we have
$$
(\sigma_{1}\Id + \mathbb{P}_{1})X(\sigma_{2}\Id +\mathbb{P}_{2})X 
 = (\sigma_{1}\Id + \mathbb{P}_{1})((1+\sigma_{2})\Id -\mathbb{P}_{1})
 = \sigma_{1}(1+\sigma_{2})\Id +(\sigma_{2} - \sigma_{1})\mathbb{P}_{1},
$$ 
and similarly
$$
(\sigma_{2}\Id + \mathbb{P}_{2})X(\sigma_{1}\Id +\mathbb{P}_{1})X = 
(\sigma_{2}\Id + \mathbb{P}_{2})((1+\sigma_{1})\Id -\mathbb{P}_{2}) = 
\sigma_{2}(1+\sigma_{1})\Id +(\sigma_{1} - \sigma_{2})\mathbb{P}_{2}.
$$ 
Using these identities for computing $\mathbb{J}_{2}^{2}$, we find
$$
(\mathbb{J}_{2})^{2} =\left[\begin{array}{cc}
\displaystyle{ \frac{\sigma_{1}}{1+\sigma_{1}}\Id + 
\frac{\sigma_{2}-\sigma_{1}}{(1+\sigma_{1})(1+\sigma_{2})}\mathbb{P}_{1} } & 0\\[5pt]
0 & \displaystyle{ \frac{\sigma_{2}}{1+\sigma_{2}}\Id + 
\frac{\sigma_{1}-\sigma_{2}}{(1+\sigma_{2})(1+\sigma_{1})}\mathbb{P}_{2}  }
\end{array}\right].
$$
The eigenvalues of $\mathbb{J}_{2}^{2}$ are thus the eigenvalues of
each of its diagonal blocks. Since the eigenvalues of the projectors
$\mathbb{P}_{j}$ are $0,1$, a direct calculation shows that the
spectrum of $(\mathbb{J}_{2})^{2}$ is $\{\sigma_{1}/(1+\sigma_{1}),
\sigma_{2}/(1+\sigma_{2})\}$, and hence we find as in the one
dimensional case
\begin{equation}\label{SpectJacobi1}
\sigma(\mathbb{J}_{2})\subset \left\{ 
+\sqrt{\frac{\sigma_{1}}{1+\sigma_{1}}}, -\sqrt{\frac{\sigma_{1}}{1+\sigma_{1}}},
+\sqrt{\frac{\sigma_{2}}{1+\sigma_{2}}}, -\sqrt{\frac{\sigma_{2}}{1+\sigma_{2}}} \right\}.
\end{equation}
Note that for $\sigma_{j}<0$ and $1+\sigma_{j}>0$, the eigenvalues
$\pm\sqrt{\sigma_{j}/(1+\sigma_{j})}$ are purely imaginary. From
(\ref{SpectJacobi1}), the spectral radius of the Jacobi method is
given by
$$
  \rho(\mathbb{J}_{2})=\max_{j=1,2}\sqrt{\Big\vert\frac{\sigma_{j}}
    {\sigma_{j}+1}}\Big\vert.
$$ 
We have thus recovered the same result as in the simple 1D model
problem from Section \ref{1dSec}. It is remarkable that the
convergence of this Jacobi iteration does neither depend on the
geometry of $\Gamma = \partial\Omega_{1} = \partial\Omega_{2}$, nor on
the dimension of the problem. {\color{black} Actually this does not
  even depend on the equation considered as all the computations
  leading to (\ref{SpectJacobi1}) are based on algebraic identities
  stemming from Proposition \ref{PropReprFormula}, which is valid at
  least for any elliptic system with piece-wise constant coefficients,
  see \cite{MR1742312} for example. Based on this observation,
  {\color{black} a} preliminary spectral analysis of multitrace
  operators for general situations {\color{black} can be found} in
  \cite{2015arXiv150800556C}. }

\begin{figure}
  \centering
  \includegraphics[width=0.9\textwidth]{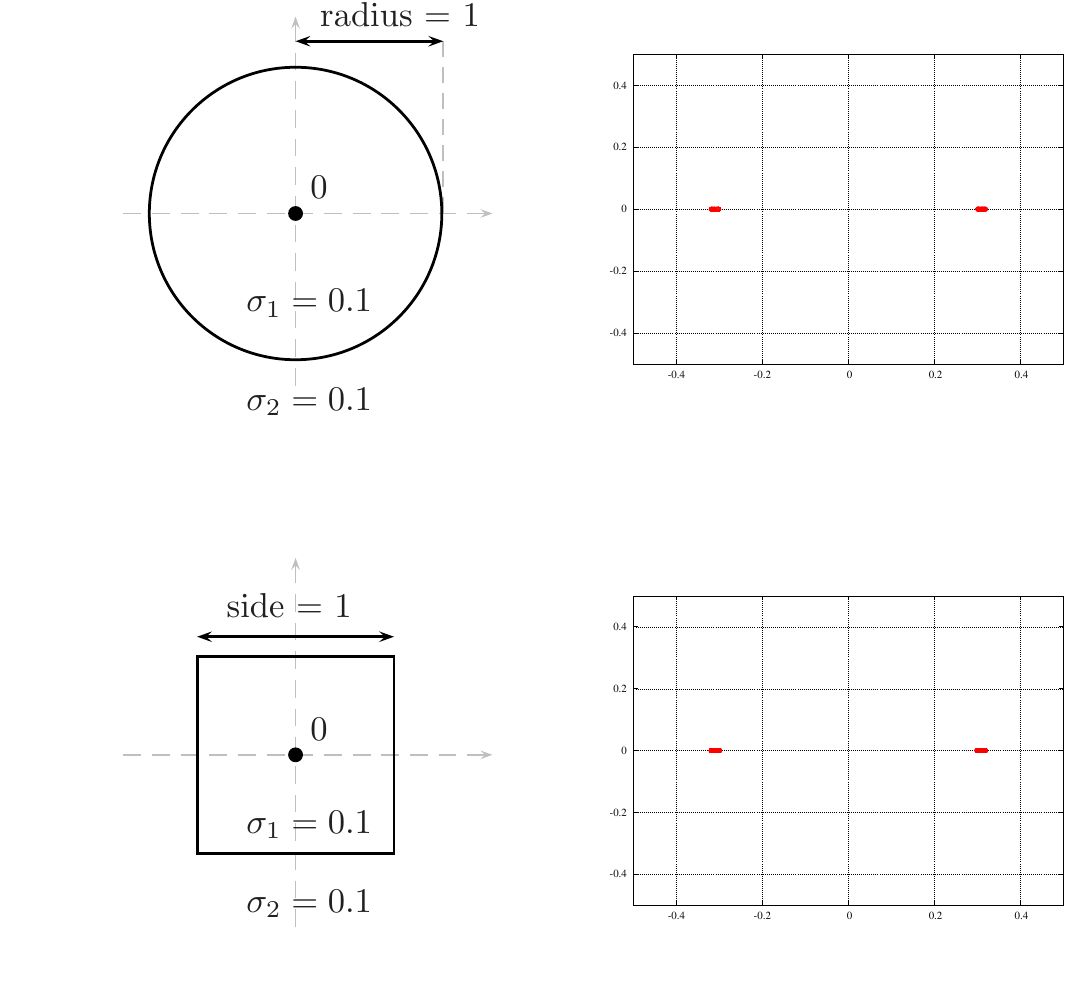}
  \caption{Spectrum in the case of the unit circle and the square: same $\sigma_j$}
\label{fig:speccircle}
\end{figure}

As an illustration, we show in Figure \ref{fig:speccircle}
a numerical approximation of the spectrum of $\mathbb{J}_{2}$ obtained
from a boundary element discretization of the Calder\'on projectors
$\mathbb{P}_{1,2}$ using $P_{1}$-Lagrange shape functions for two
different geometries: $\Gamma$ either a unit circle or a unit square.
We took $a = 1$ and $\sigma_{1} = \sigma_{2} = 0.1$ so that the exact
spectrum of $\mathbb{J}_{2}$ given in (\ref{SpectJacobi1}) is
approximately at $\{ \pm 0.301511\}$ (up to 6 digits accuracy). We
observe that the spectrum of the numerical approximation in Figure
\ref{fig:speccircle} clusters around the theoretical values $\pm
0.301511$, but is not exactly a point spectrum, which is due to
discretization and quadrature errors of the integral operators.  We
also see that the numerical spectrum appears to depend only very
weakly on the geometry.

\begin{figure}
  \centering
  \includegraphics[width=0.9\textwidth]{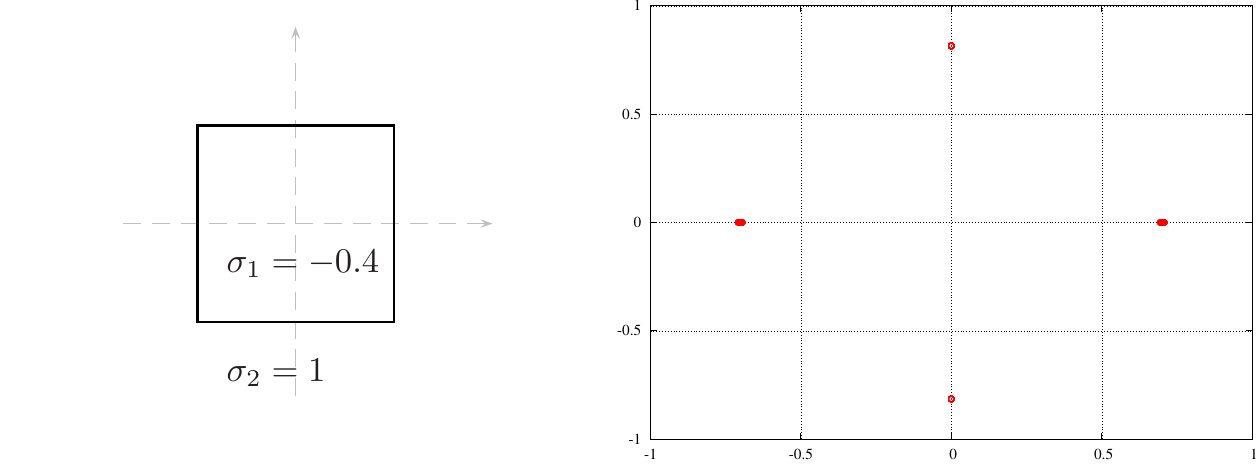}
  \caption{Spectrum in the case of the square geometry: different $\sigma_j$}
  \label{fig:specsquare}
\end{figure}

Next we consider the same computation as before with the square shaped
geometry, but in the case where the sigmas are different, $\sigma_{1}
= -0.4$ and $\sigma_{2} = 1$. We show in Figure \ref{fig:specsquare}
the corresponding spectrum of the Jacobi iteration operator.  We
clearly see that this spectrum has four clusters associated with the
two pairs of opposite eigenvalues.

\begin{figure}
\centering
\includegraphics[width=0.9\textwidth]{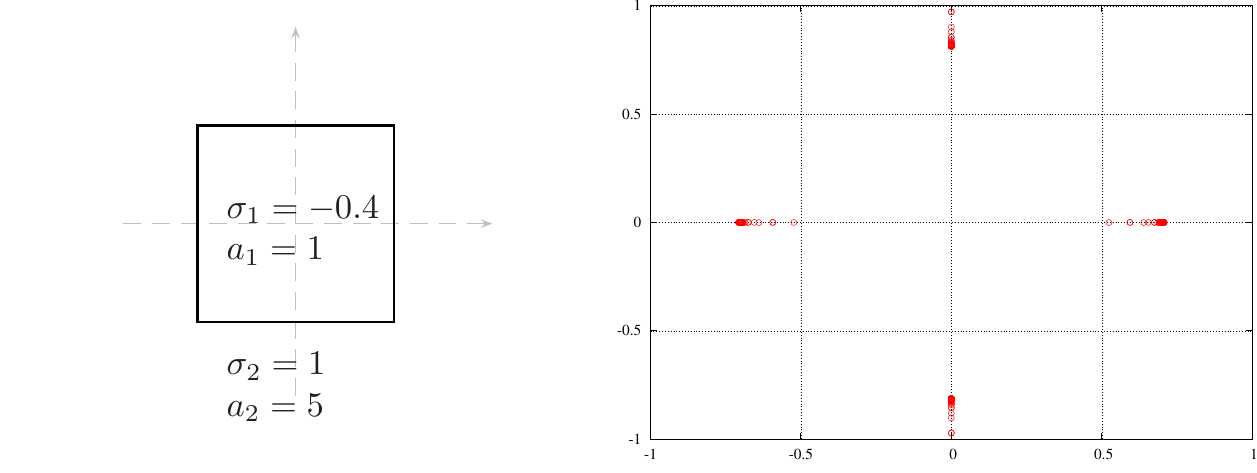}
\caption{Spectrum in the case of the square with different $\sigma_j$
  and varying coefficient $a$}
\label{fig:specheter}
\end{figure}

Finally, we consider the same experiment as above, but in the case where
the material constant $a$ is different in the two subdomains, $a_{1} =
1$ in $\Omega_{1}$ and $a_{2} = 5$ in $\Omega_{2}$. This contrast of
material characteristics only induces compact perturbations of the
integral operators, so the accumulation points of the spectrum of the
Jacobi iteration are preserved, as one can see in Figure \ref{fig:specheter}.

\subsection{Multitrace Formulation with 3 Subdomains}

\begin{figure}
\centering
\includegraphics[width=0.6\textwidth]{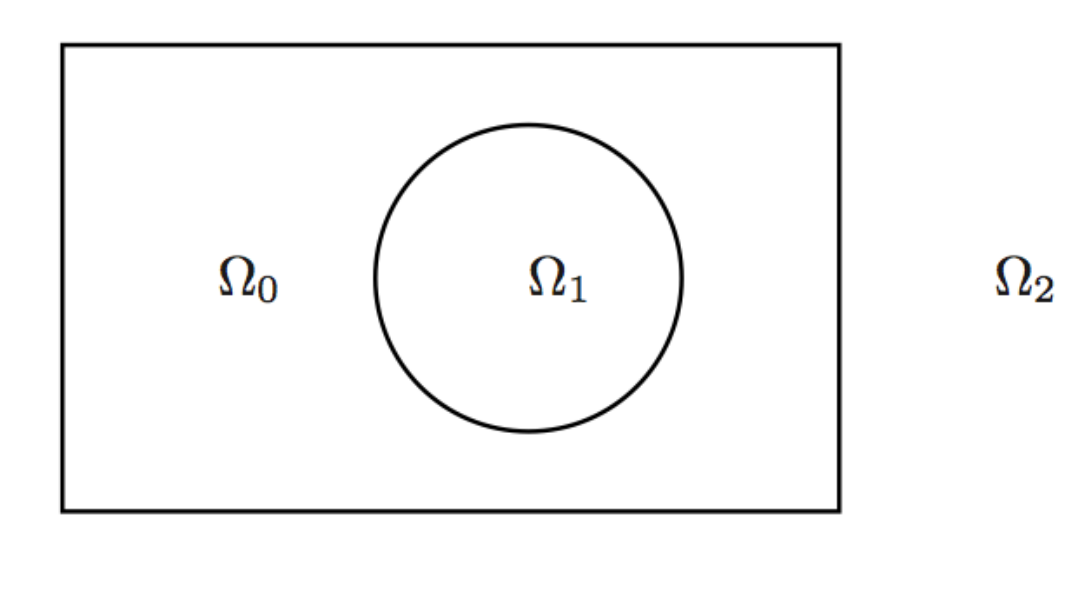}
\caption{Example of a decomposition into three subdomains}
\label{Figure3DomDec}
\end{figure}

We examine now a situation similar to Subsection
\ref{ThreeDomainDDM} in a higher dimensional context.  We consider
three Lipschitz domains with bounded boundaries $\Omega_{j},j=0,1,2$
such that $\Omega_{j}\cap \Omega_{k} = \emptyset$ for $j\neq k$ and
$\R^{d} = \overline{\Omega}_{0}\cup \overline{\Omega}_{1} \cup
\overline{\Omega}_{2}$. To fix ideas, we assume that $\Omega_{0}$ and
$\Omega_{1}$ are bounded, and $\partial\Omega_{0} =
\partial\Omega_{1}\cup \partial\Omega_{2}$ and $\partial\Omega_{1}\cap
\partial\Omega_{2} = \emptyset$, for an example, see Figure
\ref{Figure3DomDec}. 
Let $\Gamma_{1} := \partial\Omega_{0}\cap
\partial\Omega_{1}$ and $\Gamma_{2} := \partial\Omega_{0}\cap
\partial\Omega_{2}$.  For given $h^{j}_{\mathrm{D}}\in
\mH^{1/2}(\Gamma_{j})$ and $h^{j}_{\mathrm{N}}\in \mH^{-1/2}(\Gamma_{j})$,
we are interested in solving a transmission problem of the form
\begin{equation}\label{ModelEq4}
\left\{\begin{array}{l}
u\in \mH^{1}(\R^{d}\setminus( \Gamma_{1}\cup\Gamma_{2}) ),\\
-\Delta u + a^{2} u = 0\quad\textrm{in}\;\;\Omega_{j},\\
\lbr u\rbr_{\Gamma_{j}} = h_{\mathrm{D}}^{j},\quad 
\lbr \partial_{n}u\rbr_{\Gamma_{j}} = h_{\mathrm{N}}^{j},\quad\quad  j=0,1,2,
\end{array}\right.
\end{equation}
where by denoting $u_{j} = u\vert_{\Omega_{j}}$, we set $\lbr
u\rbr_{\Gamma_{j}}:= u_{j}\vert_{\Gamma_{j}} -
u_{0}\vert_{\Gamma_{j}}$ for $j=1,2$, and $\lbr
\partial_{n}u\rbr_{\Gamma_{j}}:= \bfn_{j}\cdot \nabla
u_{j}\vert_{\Gamma_{j}} + \bfn_{0}\cdot\nabla u_{0}\vert_{\Gamma_{j}}$,
$j=1,2$.

We rewrite this problem by means of a local multitrace formulation.
In what follows we denote by $T_{j}$, $j=0,1,2$ the trace operator
(\ref{DefTraceOp}) associated to each subdomain $\Omega_{j}$, and
denote the traces by $U_{j} := T_{j}(u)$.  We set $U_{0,j} :=
U_{0}\vert_{\Gamma_{j}}$. We denote by $G_{j}$ the potential
operator (\ref{DefPot}) associated to each $\Omega_{j}$, and
$\mathbb{P}_{j}$ is the corresponding Calder\'on projector.

\begin{figure}
\centering
\includegraphics[width=0.9\textwidth]{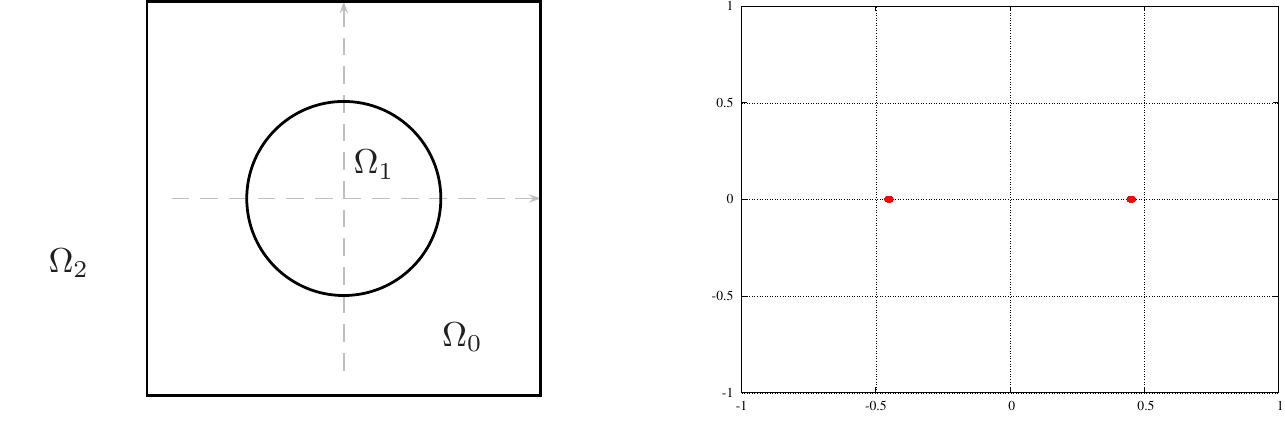}
\caption{Spectrum of the Jacobi iteration for $\sigma_{0} = \sigma_{1}
  = \sigma_{2} = 0.25$}
\label{figure6}
\end{figure}

Observe that, since we have the decomposition $\partial\Omega_{0} = 
\partial\Omega_{1}\cup \partial\Omega_{2}$, the Calder\'on projector 
$\mathbb{P}_{0}$ can be expanded into a $2\times 2$ matrix of integral 
operators, 
\begin{equation}\label{MatrixExpP0}
\mathbb{P}_{0} = 
\left\lbr\begin{array}{cc}
\widetilde{\mathbb{P}}_{1,1} & \mathbb{R}_{1,2}\\
\mathbb{R}_{2,1} & \widetilde{\mathbb{P}}_{2,2}
\end{array}\right\rbr.
\end{equation}
A close inspection of the definition of the Calder\'on projector shows
that $\mathbb{R}_{1,2} = -X\cdot (T_{1}\cdot G_{2})\cdot X$ and
similarly $\mathbb{R}_{2,1} = -X\cdot (T_{2}\cdot G_{1})\cdot X$. Take
any element $V\in \mH^{1/2}(\Gamma_{1})\times \mH^{-1/2}(\Gamma_{1})$
and set $v(\bfx):= G_{1}( X V)(\bfx)$.  This function satisfies in
particular $-\Delta v + a^{2} v = 0$ in $\Omega_{2}$, which means that
$G_{2}(T_{2}v)(\bfx) = 0$ for $\bfx\in\Omega_{1}$ according to
Proposition \ref{PropReprFormula}.  In particular $T_{1} G_{2}(T_{2}v)
= 0$. From this we conclude that
$$
\begin{array}{ll}
\mathbb{R}_{1,2}\cdot \mathbb{R}_{2,1}(V) 
& = (\, X\cdot (T_{1}\cdot G_{2})\cdot X\,)\cdot 
(\, X\cdot (T_{2}\cdot G_{1})\cdot X\,) (V)  \\[5pt]
& = X (T_{1}\cdot G_{2})\cdot  (T_{2}\cdot G_{1}) X(V)\\[5pt]
& = X \cdot T_{1}\cdot G_{2}(T_{2} v) = 0.
\end{array}
$$
We prove in a similar manner that $\mathbb{R}_{2,1}\cdot
\mathbb{R}_{1,2}(V) = 0$ for any $V\in \mH^{1/2}(\Gamma_{2})\times
\mH^{-1/2}(\Gamma_{2})$. Since, in the above arguments, $V$ was chosen
arbitrarily, we conclude that
$$
 \mathbb{R}_{2,1}\cdot \mathbb{R}_{1,2} = 0,\quad
\textrm{and} \quad
 \mathbb{R}_{1,2}\cdot \mathbb{R}_{2,1} = 0.
$$
Other remarkable identities involving $\mathbb{R}_{2,1},
\mathbb{R}_{1,2}$ can be derived.  Indeed for any $V\in
\mH^{1/2}(\Gamma_{1})\times \mH^{-1/2}(\Gamma_{1})$, the function
$v(\bfx):= G_{1}( X V)(\bfx)$ satisfies $-\Delta v + a^{2} v = 0$ in
$\Omega_{2}$, so $G_{2}(T_{2}v)(\bfx) = v(\bfx)$ for $\bfx\in
\Omega_{2}$ according to Proposition \ref{PropReprFormula} and,
consequently, $\mathbb{P}_{2}T_{2}(v) = T_{2}(v)$ which leads to
$\mathbb{P}_{2}X\mathbb{R}_{2,1}V = X\mathbb{R}_{2,1}V$. We prove
similarly $\mathbb{P}_{1}X\mathbb{R}_{1,2}V = X\mathbb{R}_{1,2}V$ for
any $V\in \mH^{1/2}(\Gamma_{2}) \times \mH^{-1/2}(\Gamma_{2})$. Since,
in this argumentation, the $V$'s were chosen arbitrarily we 
conclude that
$$
 \mathbb{P}_{1}X\mathbb{R}_{1,2} = X\mathbb{R}_{1,2},\quad
\textrm{and} \quad
 \mathbb{P}_{2}X\mathbb{R}_{2,1} = X\mathbb{R}_{2,1}.
$$
We prove in a similar manner that
$\mathbb{R}_{1,2}\widetilde{\mathbb{P}}_{2} = \mathbb{R}_{1,2}$ and
$\mathbb{R}_{2,1}\widetilde{\mathbb{P}}_{1} = \mathbb{R}_{2,1}$. For
the diagonal blocks of (\ref{MatrixExpP0}), we prove in a
similar manner as in (\ref{CommutingRelation}) that
$$
  X\mathbb{P}_{j}X = \Id - \widetilde{\mathbb{P}}_{j}.
$$
In particular the $\widetilde{\mathbb{P}}_{j}$ are projectors.  Given
three relaxation parameters $\sigma_{j}$, $j=0,1,2$, the local
multitrace formulation of Problem (\ref{ModelEq4}) is again of the
same form here as in the simple one dimensional case
(\ref{eq:multi3sd1}), and we obtain in matrix form
$$
\left\lbr\begin{array}{cccc}
(1+\sigma_{1})\Id -\mathbb{P}_{1} & -\sigma_{1}X & 0 & 0\\
-\sigma_{0}X & (1+\sigma_{0})\Id -\widetilde{\mathbb{P}}_{1} & -\mathbb{R}_{2,1} & 0\\
0 & -\mathbb{R}_{1,2} & (1+\sigma_{0})\Id -\widetilde{\mathbb{P}}_{2} & -\sigma_{0}X\\
0 & 0 & -\sigma_{2}X & (1+\sigma_{2})\Id - \mathbb{P}_{2}
\end{array}\right\rbr
\left\lbr\begin{array}{c}
U_{1}\\
U_{0,1}\\
U_{0,2}\\
U_{2}
\end{array}\right\rbr = F,
$$
where $F$ is the right hand side taking into account the 
data $h_{\mathrm{D}}^{j},h_{\mathrm{N}}^{j}$, $j=1,2$, as we have
shown in the simple 1D case. To simplify notations,
we set $\alpha_{j} := \sigma_{j}^{-1} (1+\sigma_{j})^{-1}$, so that 
$\alpha_{j}(\sigma_{j}\Id + \mathbb{P}_{j})\cdot( (1+\sigma_{j})\Id - 
\mathbb{P}_{j}) = \Id$. The Jacobi iteration matrix 
associated to the multitrace formulation then becomes
\begin{equation}\label{DefJ3}
\begin{array}{c}
\mathbb{J}_{3} = 
\left\lbr\begin{array}{cccc}
\alpha_{1}(\sigma_{1}\Id + \mathbb{P}_{1}) & 0 & 0 & 0\\
0 & \alpha_{0}(\sigma_{0}\Id + \widetilde{\mathbb{P}}_{1}) & \alpha_{0}\mathbb{R}_{1,2} & 0\\
0 & \alpha_{0}\mathbb{R}_{2,1} & \alpha_{0}(\sigma_{0}\Id + \widetilde{\mathbb{P}}_{2}) & 0\\
0 & 0 & 0 & \alpha_{2}(\sigma_{2}\Id + \mathbb{P}_{2})
\end{array}\right\rbr\\
\cdot \left\lbr\begin{array}{cccc}
0 & \sigma_{1}X & 0 & 0\\
\sigma_{0}X & 0 & 0 & 0\\
0 & 0 & 0 & \sigma_{0}X\\
0 & 0 & \sigma_{2}X & 0
\end{array}\right\rbr.
\end{array}
\end{equation}
For the sake of simplicity, to examine the spectrum of the Jacobi operator $\mathbb{J}_{3}$, 
we restrict our analysis to the case where 
$$
  \sigma_{0}  =\sigma_{1} = \sigma_{2} :=  \sigma\quad\Longrightarrow\quad
    \alpha_{0} = \alpha_{1}=\alpha_{2}=: \alpha.
$$
Under this hypothesis, we can clearly factorize $\alpha\sigma =
(1+\sigma)^{-1}$ in (\ref{DefJ3}) so it suffices to examine the
spectrum of $(1+\sigma)\mathbb{J}_{3}$. As in Subsection \ref{2Dom2D},
we study the square of this operator. Tedious, but straightforward
calculations then yield
\begin{equation}
\begin{array}{r}
(1+\sigma)^{2}(\mathbb{J}_{3})^{2} = 
\left\lbr\begin{array}{cccc}
(\sigma\Id + \mathbb{P}_{1}) & 0 & 0 & 0\\
0 & (\sigma\Id + \widetilde{\mathbb{P}}_{1}) & \mathbb{R}_{1,2} & 0\\
0 & \mathbb{R}_{2,1} & (\sigma\Id + \widetilde{\mathbb{P}}_{2}) & 0\\
0 & 0 & 0 & (\sigma\Id + \mathbb{P}_{2})
\end{array}\right\rbr\\  \cdot 
\left\lbr\begin{array}{cccc}
(1+\sigma)\Id - \mathbb{P}_{1} & 0 & 0 & X\mathbb{R}_{1,2}X\\
0 & (1+\sigma)\Id - \widetilde{\mathbb{P}}_{1} & 0 & 0\\
0 & 0 & (1+\sigma)\Id - \widetilde{\mathbb{P}}_{2} & 0\\
X\mathbb{R}_{2,1}X & 0 & 0 & (1+\sigma)\Id - \mathbb{P}_{2}
\end{array}\right\rbr\\
 = \left\lbr\begin{array}{cccc}
\sigma(1+\sigma)\Id & 0 & 0 & (1+\sigma)X\mathbb{R}_{1,2}X\\
0 & \sigma(1+\sigma)\Id & \sigma\mathbb{R}_{1,2} & 0\\
0 & \sigma\mathbb{R}_{2,1} & \sigma(1+\sigma)\Id & 0\\
(1+\sigma)X\mathbb{R}_{2,1}X & 0 & 0 & \sigma(1+\sigma)\Id
\end{array}\right\rbr.
\end{array}
\end{equation}
In the course of the above calculations, we used again several
remarkable identities: $\mathbb{P}_{2}X\mathbb{R}_{2,1} =
X\mathbb{R}_{2,1}$ and $\mathbb{P}_{1}X\mathbb{R}_{1,2} =
X\mathbb{R}_{1,2}$ as well as
$\mathbb{R}_{2,1}\widetilde{\mathbb{P}}_{1} = \mathbb{R}_{2,1}$ and
$\mathbb{R}_{1,2}\widetilde{\mathbb{P}}_{2} = \mathbb{R}_{1,2}$.  Now
observe that $(1+\sigma)(\mathbb{J}_{3})^{2} - \sigma\Id$ only
contains extra diagonal terms involving $\mathbb{R}_{1,2}$ and
$\mathbb{R}_{2,1}$. Since $\mathbb{R}_{2,1}\mathbb{R}_{1,2} =
\mathbb{R}_{1,2}\mathbb{R}_{2,1} = 0$, taking the square of this
operator yields
$$
\Big(\; (1+\sigma)\mathbb{J}_{3}^{2} -\sigma\Id\;\Big)^{2} = 
\left\lbr\begin{array}{cccc}
X\mathbb{R}_{1,2}\mathbb{R}_{2,1}X & 0 & 0 & 0\\
0 & \mathbb{R}_{1,2}\mathbb{R}_{2,1} & 0 & 0\\
0 & 0 & \mathbb{R}_{2,1}\mathbb{R}_{1,2} & 0\\
0 & 0 & 0 & X\mathbb{R}_{2,1}\mathbb{R}_{1,2}X
\end{array}\right\rbr = 0.
$$
From this we conclude that the only eigenvalue of $\mathbb{J}_{3}^{2}$
is $\sigma/(1+\sigma)$.  This gives very precise information about the
spectrum of $\mathbb{J}_{3}$ in the case where all relaxation
parameters are equal, i.e.
\begin{equation}\label{SpectreTheorique3Domaines}
  \sigma(\mathbb{J}_{3})\subset \Big\{+\sqrt{\frac{\sigma}{1+\sigma}}, 
  -\sqrt{\frac{\sigma}{1+\sigma}}\;\Big\}.
\end{equation}
Considering once again a discretization of the boundary integral 
operators by Lagrange $P_{1}$ shape functions, we show in 
Figure \ref{figure6}
the results of a numerical experiment for the geometry shown on the
left in Figure \ref{figure6}, which is a configuration with three
subdomains. We chose to solve $-\Delta u + a^{2}u = 0$ in each
subdomain with $a = 1$. We represent the spectrum of the Jacobi
iteration matrix associated to the local multitrace formulation in
Figure \ref{figure6} on the right, where all relaxation parameters in
all subdomains are equal to $\sigma = 0.25$. In accordance with
(\ref{SpectreTheorique3Domaines}), we see that eigenvalues cluster
around the pair of opposite real values $\pm \sqrt{0.25/1.25}\simeq
0.44721$.

\begin{figure}
\centering
\includegraphics[width=0.9\textwidth]{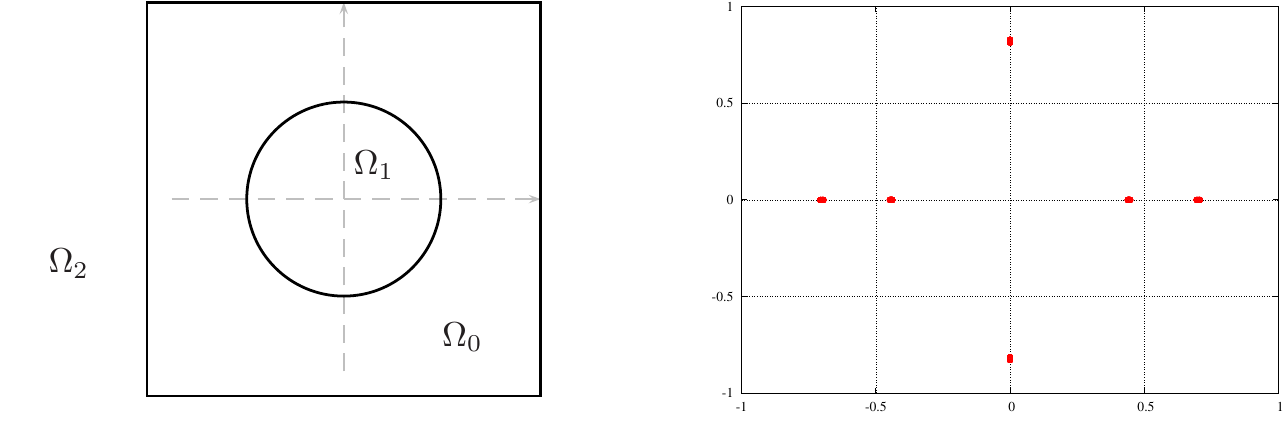}
\caption{Spectrum of the Jacobi iteration for $\sigma_{0} = -0.4$,
  $\sigma_{1} = 1$ and $\sigma_{2} = 0.25$}
\label{figure7}
\end{figure}

In Figure \ref{figure7}, 
we present the spectrum of the Jacobi iterations for a similar
numerical experiment except that we considered three different values
of the relaxation parameters, taking $\sigma_{0} = -0.4$, $\sigma_{1}
= 1$, $\sigma_{2} = 0.25$. We observe clusters of eigenvalues around
the three pairs of opposite values
$\pm\sqrt{\sigma_{j}/(1+\sigma_{j})}$ which is consistent with both
Subsection \ref{ThreeDomainDDM} and (\ref{SpectreTheorique3Domaines}).

\begin{figure}
\centering
\includegraphics[width=0.4\textwidth]{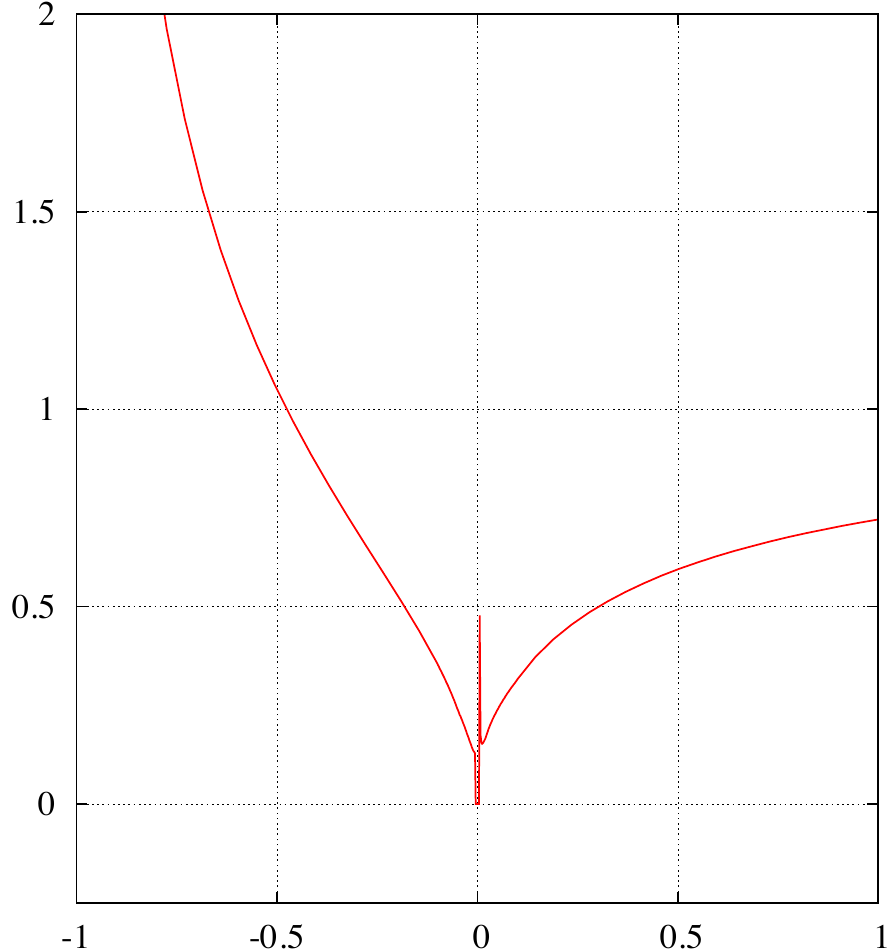}
\caption{Spectral radius of the Jacobi iteration versus $\sigma$ in
  the case where $\sigma_{0} = \sigma_{1} = \sigma_{2} = \sigma$}
\label{figure8}
\end{figure}

Finally, in Figure \ref{figure8}, 
we consider the case where all relaxation parameters are equal to
$\sigma$ and present the spectral radius of the Jacobi iteration
versus the value of this relaxation parameter $\sigma$ . We
essentially recover the curve presented in Figure \ref{fig:conv1d},
which shows that the fundamental convergence properties of the
iterative multitrace formulation we studied first on a simple one
dimensional model problem remain in this general situation. The
additional overshoot we see close to $\sigma = 0$ in Figure
\ref{figure8} compared to Figure \ref{fig:conv1d} is due to the
numerical difficulty which we explained in the one dimensional case
when $\sigma$ approaches zero.

\section{Conclusion}

We used a simple one dimensional model problem to present a recent
multitrace formulation with relaxation parameters without resorting to
a functional analysis framework. The simple setting allowed us to
study a natural block Jacobi iteration for the multitrace formulation,
and to determine the dependence of this iteration on the relaxation
parameter. We also determined an optimal choice for the relaxation
parameter, and obtained an algorithm with converges in a finite number
of steps. This algorithm is related to a well know algorithm that
also has this property: the optimal Schwarz method. We then left our
simple model problem and showed that the properties we discovered hold
also in a much more general higher dimensional setting, and this
independently of the geometry of the decomposition. An important open
question is the cost of such multitrace formulations and their
associated iterative solution. For optimal Schwarz methods it is known
that it is more efficient to use approximations of the Dirichlet to
Neumann maps to obtain practical algorithms. It is not clear yet how
in the multitrace formulation such approximations could be introduced.

\bibliographystyle{abbrv}
\bibliography{paper}
\end{document}